\documentstyle{amsart}

\newcommand{\say}[1]{\hbox{\ #1\ }}
\renewcommand{\bar}{\overline}
\renewcommand{\to}{\longrightarrow}
\renewcommand{\and}{\say{and}}

\renewcommand{\iff}{\say{if and only if}}

\newcommand{\ov}[1]{\overline{#1}}

\newcommand{\upto}[3]{#1_{#2},\dots,#1_{#3}}

\newcommand{\bbC}{{\Bbb C}}
\newcommand{\bbR}{{\Bbb R}}
\newcommand{\bbZ}{{\Bbb Z}}
\newcommand{\bbP}{{\Bbb P}}

\newcommand{\CPone}{{\Bbb C \Bbb P}^1}
\newcommand{\CPtwo}{{\Bbb C \Bbb P}^2}

\newcommand{\GL}{\hbox{\rm GL}}
\newcommand{\SO}{\hbox{\rm SO}}
\newcommand{\Orth}{\hbox{\rm O}}
\newcommand{\PSL}{\hbox{\rm PSL}}
\newcommand{\SL}{\hbox{\rm SL}}
\newcommand{\PSU}{\hbox{\rm PSU}}
\newcommand{\SU}{\hbox{\rm SU}}

\def\rank {\mathop{\hbox{\rm rank}}\nolimits}
\def\re   {\mathop{\hbox{\rm Re}}\nolimits}
\def\im   {\mathop{\hbox{\rm Im}}\nolimits}
\def\supp {\mathop{\hbox{\rm supp}}\nolimits}
\def\Spin {\mathop{\hbox{\rm Spin}}\nolimits}
\def\taut {\mathop{\hbox{\rm Taut}}\nolimits}
\def\Arf  {\mathop{\hbox{\rm Arf}}\nolimits}
\def\trace{\mathop{\hbox{\rm tr}}\nolimits}
\def\res  {\mathop{\hbox{\rm res}}\nolimits}
\def\qres {\mathop{\hbox{\rm qres}}\nolimits}
\def\ord  {\mathop{\hbox{\rm ord}}\nolimits}

\def\conj {\mathop{\hbox{\rm Conj}}\nolimits}

\def\calF{{\cal F}}
\def\calH{{\cal H}}

\def\calK{{\cal K}}
\def\calO{{\cal O}}
\def\calM{{\cal M}}

\let\om =\omega
\let\sig=\sigma
\let\al =\alpha
\let\gam=\gamma
\let\ep =\varepsilon

\newcommand{\suchthat}{\;|\;}
\newcommand{\iffspace}{\;\;\;\iff\;\;\;}

\newcommand{\andspace}{\;\;\;\hbox{and}\;\;\;}
\newcommand{\spaceout}{\hbox{\ \ \ }}

\newcommand{\skewform}{\Omega}
\newcommand{\trans}[1]{#1^t}
\newcommand{\disp}{\displaystyle}

\newcommand{\arraystrut}{\rule[-3ex]{0ex}{1ex}}
\newcommand{\arraystruti}{\rule{0ex}{3ex}}
\newcommand{\boxwidth}{40ex}

\newcommand{\image}{\hbox{\,image\ }}
\newcommand{\genus}{\hbox{\,genus\,}}

\newcommand{\pfaffian}{\hbox{\,pfaffian\,}}
\renewcommand{\to}{\longrightarrow}

\newcommand{\skipline}{\vspace{.1in}}

\newcommand{\tensor}{\otimes}
\newcommand{\Moebius}{M\"obius\ }
\newcommand{\andeq}{&=&}
\newcommand{\andin}{&\!\in\!\!\!&}
\newcommand{\Gr}{\hbox{\,\rm Gr\,}} 
\newcommand{\ooo}{\hbox{\rm o}}
\newcommand{\partialslash}{{\partial\!\!\!/}}
\newcommand{\partialbar}{{\bar{\partial}}}

\newtheorem{theorem}{Theorem}
\newtheorem{lemma}{Lemma}

\newtheorem{definition}{Definition}

\newtheorem{lemmaU}{Lemma}

\newtheorem{definitionU}{Definition}

\begin{document}

%
%

%
%

\typeout{_______________________________________________ title.tex}

\pagenumbering{roman}

\title{
The Spinor Representation of Surfaces in Space}

\author{
Rob Kusner and Nick Schmitt}

\date{}

\maketitle

\begin{abstract}

The spinor representation is developed for conformal immersions of
Riemann surfaces into space.  We adapt the approach of Dennis Sullivan
\cite{Sullivan}, which treats a spin structure on a Riemann surface
$M$ as a complex line bundle $S$ whose square is the canonical line
bundle $K=T(M)$.  Given a conformal immersion of $M$ into $\bbR^3$,
the unique spin strucure on $S^2$ pulls back via the Gauss map to a
spin structure $S$ on $M$, and gives rise to a pair of smooth sections
$(s_1,s_2)$ of $S$.  Conversely, any pair of sections of $S$
generates a (possibly periodic) conformal immersion of $M$ under a
suitable integrability condition, which for a minimal surface is
simply that the spinor sections are meromorphic.  

A spin structure $S$ also determines (and is determined by) the
regular homotopy class of the immersion by way of a $\bbZ_2$-quadratic
form $q_S$.  We present an analytic expression for the Arf invariant
of $q_S$, which decides whether or not the correponding
immersion can be deformed to an embedding.  The Arf invariant also
turns out to be an obstruction, for example, to the existence of
certain complete minimal immersions.

The later parts of this paper use the spinor representation to
investigate minimal surfaces with embedded planar ends.  In general,
we show for a spin structure $S$ on a compact Riemann surface $M$ with
punctures at $P$ that the space of all such (possibly periodic) minimal
immersions of $M\setminus P$ into $\bbR^3$ (upto homothety) is the the
product of $S^1\times H^3$ with the Grassmanian of 2-planes in a
complex vector space $\calK$ of meromorphic sections of $S$.  An
important tool -- a skew-symmetric form $\Omega$ defined by residues
of a certain meromorphic quadratic differential on $M$ -- lets us
compute how $\calK$ varies as $M$ and $P$ are varied.  Then we apply
this to determine the moduli spaces of planar-ended minimal spheres
and real projective planes, and also to construct a new family of
minimal tori and a minimal Klein bottle with 4 ends.  These surfaces
compactify in $S^3$ to yield surfaces critical for the
\Moebius invariant squared mean curvature functional $W$.  On the
other hand, Robert Bryant \cite{Bryant1} has shown all $W\!$-critical
spheres and real projective planes arise this way.  Thus we find at
the same time the moduli spaces of $W\!$-critical spheres and real
projective planes via the spinor representation.
\end{abstract}

\vspace{.5 in}

Department of Mathematics

Center for Geometry, Analysis, Numerics and Graphics

Univerity of Massachusetts, Amherst, MA 01003

\vfill\pagebreak

\mbox{Our work at GANG was supported in part by NSF grants 
DMS 93-12087 and 94-04278. }
\vfill\pagebreak 
\tableofcontents

%
%

%
%

\typeout{_______________________________________________ intro.tex}

\pagenumbering{arabic}
\section*{Introduction}

In this paper we investigate the interplay between spin structures
on a Riemann surface $M$ and immersions of $M$ into three-space.
Here, a spin structure is a complex line bundle $S$ over $M$ such
that $S\tensor S$ is the holomorphic (co)tangent bundle $T(M)$ of
$M$. Thus we may view a section of a $S$ as a ``square root'' of a
$(1,0)$-form on $M$. Using this notion of spin structure, in the
first part of this paper we develop the notion of the {\em spinor
representation of a surface in space}, generalizing an observation
of Dennis Sullivan \cite{Sullivan}. The classical Weierstrass
representation for a minimal surface is
\[
(g,\eta)\to \re\int(1-g^2,i(1+g^2),2g)\eta,
\]
where $g$ and $\eta$ are respectively a meromorphic function and
1-form on the underlying compact Riemann surface.  The spinor
representation (Theorem \ref{induced spin structure}) is
\[
(s_1,s_2)\to\re\int(s_1^2-s_2^2,i(s_1^2+s_2^2),2 s_1 s_2),
\]
where $s_1$ and $s_2$ are meromorphic sections of a spin structure
$S$. Either representation gives a (weakly) conformal harmonic map
$M\rightarrow\bbR^3$, which therefore parametrizes a (branched)
minimal surface. In fact, either can be used to construct {\em any}
conformal immersion (not necessarily minimal) of a surface if we
relax the meromorphic condition on the data to a suitable
integrability condition: in terms of spinors, we require that
$(s_1,\,s_2)$ satisfy the first-order equation
\[
\partialbar(s_1,\,s_2) = H(|s_1|^2+|s_2|^2)(\bar{s}_2,\,-\bar{s}_1),
\]
which clearly reduces to the Cauchy-Riemann equations on a minimal surface
(Theorem \ref{dirac}), that is, when the mean curvature $H$ vanishes.

One feature of the spinor representation is that fundamental
topological information, such as the regular homotopy class of the
immersion, can be read off directly from the analytic data (Theorem
\ref{homotopy theorem}).
In fact for tori, and more generally for 
hyperelliptic Riemann surfaces (Theorem \ref{hyperelliptic}), we are able to give an explicit calculation of the Arf
invariant for the immersion: the Arf
invariant
distinguishes whether or not an immersion of an orientable surface
is regularly homotopic to an embedding.  We also consider
the spinor representation for nonorientable surfaces in terms
of a lifting to the orientation double cover
(Theorem \ref{nonorientable theorem}).  This is
sufficient for constructing minimal examples later in the paper, but is less
satisfying theoretically. In a future paper, we plan to consider the
general case from the perspective of ``pin'' structures, and also give a
more direct differential geometric treatment of the Arf invariant.

The second part of this paper considers finite-total-curvature minimal
surfaces from the viewpoint of the spinor
representation, particularly surfaces with embedded planar ends.
It is well-known (see \cite{Bryant1},
\cite{Kusner1}, \cite{Kusner2}) that such surfaces conformally
compactify to give extrema for the squared mean curvature integral
$W=\int H^2dA$ popularized by Willmore.  Conversely, for genus zero,
all $W\!$-critical surfaces arise this way \cite{Bryant1}.  Using the
spinor representation to study these special minimal surfaces has the
computational advantage of converting certain quadratic conditions to
linear ones.  In fact, associated to a spin structure $S$ on a closed
orientable Riemann surface $M$ is a vector space $\calK$ of sections
of $S$ such that pairs of independent sections $(s_1,s_2)$ from
$\calK$ form the spinor representations of all the minimal immersions
of $M$ with embedded planar ends (Theorems \ref{K condition} and \ref{theorem: moduli space}).  Thus
the problem of finding all these immersions is reduced to an algebraic
problem (Theorem \ref{p:space}); to better understand $\calK$, a
skew-symmetric bilinear form $\skewform$ is defined from whose kernel
$\calK$ is computable (Definition \ref{skewform} and Theorem 12).

The third (and final) part of this paper is devoted to the
construction of examples and to classification results.  Specifically,
for a given finite topological type of surface, we explore the moduli
space $\calM$ of immersed minimal surfaces (up to similarity) of this
type with embedded planar ends: the dimension and topology of $\calM$,
convergence to degenerate cases (that is, the natural closure of
$\calM$), and examples with special symmetry (which correspond to
singular points of $\calM$).  The tools mentioned above permit the
broad outline of a solution, but require ingenuity to apply in
particular cases.  For example, the form $\skewform$ allows the moduli
space to be expressed as a determinantal variety which determines how
the location of the ends can vary along the Riemann surface $M$.
However, this determinantal variety is only computable when the number
of ends is small.  Furthermore, the basic tools, being algebraic
geometrical, ignore the real analytic problems of removing periods and
branch points.  The latter require much subtler and often {\em ad hoc}
methods.

Previously known results concerning minimal surfaces
with embedded planar ends include the following:
\begin{enumerate}

\item[$\bullet$]
spheres exist for $4$, $6$, and every $n\ge 8$ ends \cite{Bryant1},
\cite{Kusner2}, \cite{Peng};

\item[$\bullet$]
there are no immersed spheres with $3$, $5$, and $7$ ends \cite{Bryant2};

\item[$\bullet$]
the moduli spaces of immersed spheres with $4$ and $6$ ends, and
projective planes with 3 ends have been determined \cite{Bryant2};

\item[$\bullet$]
there exist rectangular tori with 4 ends \cite{Costa}.
\end{enumerate}

\noindent Using the spinor representation we find:
\begin{enumerate}

\item[$\bullet$] a new proof of the non-existence of spheres with $3$,
$5$ and $7$ ends is given using the skew-symmetric form $\skewform$
(Theorem \ref{theorem:nogenus0});

\item[$\bullet$] the moduli space of spheres with $2p$ ends
($2\le p\le 7$) is shown to be $4(p-1)$-dimensional (Theorem
\ref{theorem:small ends});

\item[$\bullet$]
the point which compactifies the moduli space of projective planes
with 3 ends is proved to be a \Moebius strip, and all symmetries of
these surfaces are found (Theorem \ref{theorem:projective});

\item[$\bullet$]
there are no three-ended tori (Theorem \ref{theorem:3ends}); 

\item[$\bullet$]
there is a real two-dimensional family of four-ended immersed examples
on each conformal type of torus (Theorem \ref{theorem:4ends2});

\item[$\bullet$]
there exists an immersed Klein bottle with four ends 
(Theorem \ref{theorem:klein}).
\end{enumerate}

For higher genus, the general methods we have developed here also yield
(possibly branched) minimal immersions with embedded planar ends, but
it becomes more and more difficult to determine precisely when branch
points are absent or periods vanish: we again postpone this case to a
future paper.

Most of the theorems presented here were worked out while we visited
the Institute for Advanced Study during the 1992 Fall term, and were
first recorded in \cite{Schmitt}.  We thank the
School of Mathematics at the Institute for its hospitality, as well as
A. Bobenko, G. Kamberov, P. Norman, F. Pedit, U.  Pinkall, J.
Richter, D. Sullivan, J. Sullivan and I. Taimanov for
their comments and interest. Additionally, we should mention some
more recent related developments in \cite{Bobenko} and
\cite{KamberovGANG}.

\vfill\pagebreak

%
%

\typeout{_______________________________________________ part1.tex}

\part{Spinors, Regular Homotopy Classes and the Arf Invariant}

The notion of a spin structure is developed and used to describe the
spinor representation of a surface in space. Section \ref{define q}
defines a ``quadratic form'' which can be used to completely classify
the spin structures on a surface, and Section \ref{2sphere} computes
coordinates for the unique spin structure on the Riemann sphere. In
the next two sections, the spinor representation of a surface is
explained and related to the regular homotopy class of the surface.
Section \ref{holo} shows equivalent characterizations of spin
structures, the most useful of which will be that of representing spin
structures by holomorphic differentials. These differentials are
computed on tori (and, in Appendix B, on hyperelliptic Riemann surfaces).
Section \ref{group action} takes up the question of group action
on spinors, and computes
the group which performs Euclidean similarity transformations. Two
surfaces which are transforms of each other under the action of this
group are considered to be the same. The final two sections discuss
briefly the technicalities of periods and nonorientable surfaces.

\section{Spin structures on a surface}\label{spin structures}

A spin structure on an $n$-dimensional (spin) manifold $M$ is a
certain two-sheeted covering map of the $\SO(n)$-frame-bundle on $M$ to a
$\Spin(n)$-bundle (see \cite{Milnor}, \cite{Lawson}). When $ n=2 $, this notion of
spin structure may easily be reduced to the following definition in
terms of a quadratic map between complex line bundles:
\vspace{-0.15in}

\begin{figure}[htbp]
\setlength{\unitlength}{1.8cm}

\centering\begin{picture}(2,1.4)(-.5,0)
\small
\put(0.2,0.95){\vector(1,0){0.6}}
\put(1.1,0.8){\vector(0,-1){0.6}}
\put(0.2,0.8){\vector(1,-1){0.7}}

\put(.95,-.05){$M$}
\put(-0.1,0.9){$S$}
\put(0.9,0.9){$K = T(M)$}
\put(.4,1.1){$\mu$}

\end{picture}
\caption{Spin structure}
\end{figure}

\vspace{-0.1in}
\begin{definition}\label{define spin structure}
A {\em spin structure} on a Riemann surface $M$ is a complex line bundle
$S$ over $M$ together with a smooth surjective fiber-preserving map 
$\mu:S\to K$ to the holomorphic (co)tangent bundle $K = T(M)$ satisfying
\begin{equation}\label{spin condition}
\mu(\lambda s) = \lambda^2\mu(s)
\end{equation}
for any section $s$ of $S$. 
We refer to a section of $S$ as a {\em spinor}.
\end{definition}

Two spin structures $(S,\mu)$ and $(S',\mu')$ on a Riemann surface $M$
are {\em isomorphic} if there is a line bundle isomorphism
$\delta:S\to S'$ for which $\mu = \mu' \delta$. Hence two spin
structures may be isomorphic as line bundles and yet not be
isomorphic as spin structures. The number of non-isomorphic spin
structures on a Riemann surface $M$ is equal to the cardinality of
$H^1(M,\bbZ_2)$. (This count remains true for spin manifolds in general:
see \cite{Milnor}.)  In particular, if $M$ is a closed Riemann surface
of genus $g$, there are $2^{2g} = \# H^1(M,\bbZ_2)$ such structures on
$M$.

An important example is the annulus $A = \bbC^*$.
There are exactly two non-isomorphic spin structures on $A$,
which can be given explicitly as follows.
The (co)tangent bundle $T(A)$ may be identified with 
$A\times\bbC$ by means of the global trivialization
\[
 a \left.dz\right|_p \mapsto (p,a).
\]
Let $S_0=S_1=A\times\bbC$ and define maps $\mu_k:S_k\to T(A)$ for $k=0,1$ by
\[
\begin{array}{lcl}
\mu_0(z,w)=(z,w^2),\\
\mu_1(z,w)=(z,zw^2).
\end{array}
\]
Then $(S_k,\mu_k)$ are spin structures on $A$
since $\mu_k$ satisfies the condition (\ref{spin condition}).
Though $S_0$ and $S_1$ are isomorphic line bundles over $A$, they are
non-isomorphic spin structures. For if $S_0$ and $S_1$ were
isomorphic spin structures with bundle isomorphism 
$\delta:S_0\to S_1$ satisfying $\mu_0 = \mu_1\delta$, then $\delta$
would be of the form $(z,w)\mapsto (z,f(z,w))$. Then
$w^2 = z f^2$, implying that
$z$ has a consistent square root on $\bbC^*$, which is impossible.

\section{The quadratic form associated to a spin structure}\label{define q}

In this section, the Riemann surface $M$, its holomorphic (co)tangent bundle, and
the spin structure are replaced with the corresponding real
manifold and real vector bundles. In particular,  all vector fields in this
section are {\em real} vector fields.

To each spin structure $S$ on the Riemann surface $M$ we associate a 
$\bbZ_2$-valued quadratic form
\[
q_S:H_1(M,\bbZ_2)\to \bbZ_2.
\]
To say that $q_S$ is quadratic means that for all $c_1$, $c_2\in H_1(M,\bbZ_2)$
we have
\[
q_S(c_1 + c_2) = q_S(c_1) + q_S(c_2) + c_1\cdot c_2.
\]
where $c_1\cdot c_2$ denotes the intersection number (mod 2) of
$c_1$ with $c_2$.

To define $q_S(c)$, let $\al:S^1\to M$ be a smooth embedded representative of $c$
(the existence of such an $\al$ follows from results in \cite{Meeks2}).
Let $v$ be a smooth vector field along $\al$ which lifts to a section of $S$
along $\al$, and let $w(\al,\,v)$ denote the total turning number
(mod 2) of the derivative vector $\al'$ against $v$ along $\al$.
Define $q_S(c) = w(\al,\,v)+1$.

\begin{theorem}\label{quadratic theorem}
The form $q_S:H_1(M,\bbZ_2)\to \bbZ_2$ is
well-defined, that is, independent of the choice of the vector field $v$
and the choice of embedded representative $\al$, and $q_S$ is quadratic
in the above sense.
\end{theorem}

The proof is given in Appendix \ref{quadratic appendix} (see also
\cite{Atiyah} and \cite{Johnson}).

A well-known result (see, for example, \cite{Pinkall})
is that the equivalence class of a quadratic form 
$q:H_1(M,\bbZ_2)\to\bbZ_2$ under linear changes of bases of $H_1(M,\bbZ_2)$
is determined by its {\em Arf invariant}
\begin{equation}\label{Arf}
\Arf q = \frac{1}{\sqrt{\# H}}\sum_{ \al\in H} (-1)^{q(\al)},
\end{equation}
where $H = H_1(M,\bbZ_2)$.
The quadraticity of $q$ insures that
this invariant has values in $\{+1,-1\}$.
For a compact surface of genus $g$,
there are $2^{2g-1}+2^{g-1}$ spin structures for which the 
$\Arf$ invariant of the corresponding quadratic form is $+1$, and
$2^{2g-1}-2^{g-1}$ spin structures for which it is $-1$
(compare Appendix \ref{hyper appendix}).

An alternate interpretation \cite{Arbarello} of the Arf invariant 
is $\Arf q_S = (-1)^{\dim \calH_S},$
where $\calH_S$ is the space of holomorphic sections of $S$ (also see
\cite{Atiyah}).

\section{The spin structure on the Riemann sphere}\label{2sphere}

The following description of the unique spin structure on $S^2$, as
well as the spinor representation of a surface in the next section, are
adapted from \cite{Sullivan}.
Identify
\[
S^2 \cong [Q] = \{[z_1,z_2,z_3]\in\CPtwo\suchthat
     z_1^2+z_2^2+z_3^2=0\},
\]
where $Q$ is the null quadric
\[
Q = \{(z_1,z_2,z_3)\in\bbC^3\suchthat z_1^2+z_2^2+z_3^2=0\}.
\]
Then $T(S^2)$ may be identified with the restriction to $[Q]$ of the
tautological line bundle
\[
\taut(\CPtwo) = \{(\Lambda,x)\in\CPtwo\times\bbC^3\suchthat
x\in\Lambda\}
\]
(here, $\CPtwo$ is thought of as the lines in $\bbC^3$), so
\begin{equation}\label{T(S^2)}
T(S^2) \cong \taut(\CPtwo)|_{[Q]} = \{(\Lambda,x)\in[Q]\times Q\suchthat
     x=0\hbox{ or } \pi(x)\in\Lambda\},
\end{equation}
where $\pi:Q\to[Q]$ is the canonical projection.
Given this, the unique spin structure $\Spin(S^2)$ on $S^2$ may then
be identified with the tautological line bundle
\begin{equation}\label{Spin(S^2)}
\Spin(S^2) \cong \taut(\CPone) \cong
     \{(\Lambda,x)\in\CPone\times\bbC^2\suchthat x\in\Lambda\},
\end{equation}
with the associated mapping $\mu$ given by
\[
\mu( [z_1,z_2], (s_1,s_2) ) = ([\sig(z_1,z_2)],\sig(s_1,s_2)),
\]
where $\sig:\bbC^2\to Q$ is the ``Segre'' map defined by
\begin{equation}\label{sigma}
\sig(z_1,z_2) = (z_1^2-z_2^2,i(z_1^2+z_2^2),2z_1z_2).
\end{equation}
As may be checked, the map $\mu$ satisfies the conditions of Definition
\ref{define spin structure}.

When $T(S^2)$ and $\Spin(S^2)$ are restricted respectively to their
nonzero vectors and nonzero spinors, they have single coordinate charts
\[
\begin{array}{rcl}
\left\{\hbox{nonzero vectors in $T(S^2)$}\right\} &\to& Q\setminus\{0\}\\
\left\{\hbox{nonzero spinors in $\Spin(S^2)$}\right\} &\to& 
     \bbC^2\setminus\{0\}
\end{array}
\]
defined by taking the second component in each of (\ref{T(S^2)}) and
(\ref{Spin(S^2)}) respectively.
In this case, $\mu$ may be thought of as the two-to-one covering map
$\sig:\bbC^2\setminus\{0\}\to Q\setminus\{0\}$.

\section{The spinor representation of a surface in space}\label{spinor representation}

To describe the spinor representation, let $M$ be a connected Riemann
surface with a local complex coordinate $z$, and $X:M\to\bbR^3$ a conformal
(but not necessarily minimal) immersion of $M$ into space. Since $X$
is conformal, its $z$-derivative $\partial X = \om$ can be viewed as
a null vector in $\bbC^3$, or via (\ref{T(S^2)}), as a map into the
(co)tangent bundle $T(S^2)$. The Gauss map $g$ associated to $X$ 
can be viewed as a (not necessarily meromorphic) function
$g:M \to \bbC\cup\{\infty\}$ by identifying $S^2$ and
$\bbC\cup\{\infty\}$ (via stereographic projection).
This induces the bundle map $(\om,g)$ as in the lower square of Figure 2.

The {\em Weierstrass representation} of the immersion $X$ above is the
pair $(g,\eta)$, where $g$ is the stereographic projection of the
Gauss map, and $\eta$ is the $(1,0)$-form (again, not necessarily
meromorphic) on $M$ satisfying
\[
\partial X = \om = ( 1-g^2, i(1+g^2), 2g )\,\eta.
\]

Reversing this procedure (up to the problem of periods --- see
Section \ref{periods}) one obtains the following classical result.
\begin{theorem}
Given a bundle map $(\om,g)$ of $K = T(M)$ into $T(S^2)$,
if $\om$ satisfies the integrability condition
\[
\makebox{\em Re\,}d\om = 0,
\]
the $\bbR^3$-valued form $\hbox{\em Re\,}\om$ is closed (so locally
exact), and thus
\[
X = \makebox{\em Re\,}\int\om:M\to\bbR^3
\]
is a (possibly periodic, branched) conformal immersion with Gauss map $g$.
\end{theorem}

\vspace{-0.2in}
\begin{figure}[htbp]
\setlength{\unitlength}{1.8cm}
\centering\begin{picture}(2,2.3)(-.75,0)
\small
\multiput(0.1,-0)(0,1){3}{\vector(1,0){0.8}}
\multiput(-.3,0.8)(1.4,0){2}{\vector(0,-1){0.6}}
\multiput(-.3,1.8)(1.4,0){2}{\vector(0,-1){0.6}}

\put(-0.4,-.08){$M$}
\put(.95,-.08){$S^2$}
\put(-0.9,.95){$K=T(M)$}
\put(1,.95){$T(S^2)$}
\put(-0.35,2){$S$}
\put(1,2){$\Spin(S^2)$}
\put(.4,0.1){$g$}
\put(.4,1.1){$\om$}
\put(.4,2.1){$\psi$}
\put(-0.6,1.5){$\mu$}
\put(1.2,1.5){$\sig$}

\end{picture}
\caption{Spinor representation of a surface}
\end{figure}

The {\em spinor representation} of the immersion is obtained by lifting
(see Figure~2) $\om$ to the spin structures on $M$ and $S^2$.

\begin{theorem}\label{induced spin structure}
Let $S$ be a spin structure on $M$, and $(\psi,\,g)$ a bundle map as in
Figure~$2$. Assume the integrability condition.
Then there exists a (possibly periodic, branched)
immersion $X:M\to \bbR^3$ with Gauss map $g$,
whose differential $\partial X = \om$ lifts to $\psi$.

On the other hand, if $(\om,\,g)$ is a bundle map of $K = T(M)$ into
$T(S^2)$, then
\begin{enumerate}
\item there is a unique spin structure $S$ on $M$ such that $\om$ lifts to
a bundle map \em $\psi:S\to \makebox{\rm Spin}(S^2)$\em ;
\item there are exactly two such lifts $\psi$, and these differ only
by sign.
\end{enumerate}
\end{theorem}

\begin{pf}
The integrability condition $\re d\om = 0$ (or its spinor equivalent, in
Theorem \ref{dirac} below) implies that $X$ is well-defined up to periods.

{\it (i).}
Considering $\Spin(S^2)$ as a $\bbZ_2$-bundle on
$T(S^2)$ when restricted to nonzero spinors and vectors
respectively, let $S$ be the (unique) pullback bundle of $\Spin(S^2)$
under $\om$, and $\mu$, $\psi$ as shown. Extend $S$, $\psi$, and $\mu$
to include the zero spinors.

{\it (ii).}  
If $\iota:\Spin(S^2)\to \Spin(S^2)$ is the order-two deck transformation
for the covering $\Spin(S^2)\to T(S^2)$, then $\iota\circ\psi$ is
another map which in place of $\psi$ makes the diagram commute.
Conversely, if $\zeta:S\to \Spin(S^2)$ is such a map, then for $x\in S$,
$\zeta(x)$ is $\psi(x)$ or $\iota\circ\psi(x)$ and continuity
implies that $\zeta = \psi$ or $\iota\psi$.
\end{pf}

The Weierstrass and spinor representations are related by the
equation
\[
\om = \sig(\psi) = (s_1^2 - s_2^2, i(s_1^2+s_2^2), 2 s_1 s_2),
\]
where $\psi = (s_1,\,s_2)$ is viewed as a pair of sections of $S$, 
and the squaring-map $\mu$ is kept implicit by writing $s^2$ for
$\mu(s)$ and $st$ for $\frac{1}{4}(\mu(s+t) - \mu(s-t))$.
Thus these representations satisfy
\[
\eta = s_1^2\andspace g = s_2/s_1.
\]

How is the integrability condition expressed using spinors?  We compute in
the local coordinate $z$, where there are two sections of $S$ whose
images under $\mu$ are the (1,0)-form $dz$. Choose one
of these sections, and refer to it consistently as
$\sqrt{dz}=\varphi.$  Then any spinor can be written locally in the
form $s = f \varphi$, with $s^2 = \mu(s) = f^2 dz$.  We define
\[
\partial s = \partial \! f \,\varphi \andspace
\partialbar s = \partialbar \! f \, \varphi,
\]
sections of $K\tensor S\cong S\tensor S\tensor S$ and 
$\bar{K}\tensor S\cong \bar{S}\tensor \bar{S}\tensor S =
\bar{S}\tensor |S|^2$, respectively.
For the spinor pair $\psi = (s_1, s_2)$, we also write
$\partialslash\psi = \partial(-\bar{s}_2,\,\bar{s}_1)$, as suggested by
Kamberov \cite{KamberovGANG}; upto a conformal factor, $\partialslash$
is the {\it Dirac operator} associated to $S$.

The {\em first fundamental form}
(or equivalently, the {\em area form}) of the conformal immersion $X$ is a section of $|K|^2 \cong |S|^4$ given by
\[
2|\om|^2 = 4(|s_1|^2 + |s_2|^2)^2 = 4|\psi|^4.
\]
The {\em second fundamental form} has trace-free part
(the {\em Hopf differential}, a section of $K\tensor K$) given by
\[
\Phi = \eta\,\partial g = s_1 \partial s_2 - s_2 \partial s_1 =
-\psi\cdot\partialslash\bar{\psi}.
\]
Half its trace (the {\em mean curvature}) is
\[
H = \frac{n\cdot\partialbar\om}{ |\om|^2 } =
    \frac{- s_1 \partialbar s_2 + s_2 \partialbar s_1}{(|s_1|^2 + |s_2|^2)^2} =
    \frac{\psi\cdot\bar{\partialslash\psi}}{|\psi|^4}
\]
where now the Gauss map is viewed as the unit normal vector to the
surface
\[
n = \frac{\left(2g,\,|g|^2 - 1\right)}{|g|^2 + 1} =\frac{\left(s_1\bar{s}_2 + s_2\bar{s}_1,\,i(s_1\bar{s}_2 - s_2\bar{s}_1),\,|s_2|^2 - |s_1|^2\right)}{|s_1|^2 + |s_2|^2} \in \bbC\times\bbR = \bbR^3.
\]
Differentiating the relation $\om=\sigma(\psi)$ and using the formulas
above allow us to re- express the integrability condition {$\re d\om =
0$} as follows.

\begin{theorem}\label{dirac}
The integrability condition for the spinor representation $\psi=(s_1,\,s_2)$
is that the matrix
\[
\left(
\begin{array}{cccc}
s_1 & s_2 & \bar{s}_1 & \bar{s}_2 \\
-\partial\bar{s}_2 & \partial\bar{s}_1 &
-\partialbar s_2 & \partialbar s_1
\end{array}
\right)
=
\left(
\begin{array}{cc}
\psi & \bar{\psi}\\
\partialslash\psi & \bar{\partialslash\psi}
\end{array}
\right)
\]
has (real) rank 1.
Equivalently, $\psi$ satisfies a non-linear Dirac equation
\[
\partialslash\psi = H |\psi|^2 \psi
\]
for some real-valued function $H$, necessarily the mean curvature of
the surface.
\end{theorem}
\noindent
Versions of this last equation have been observed also by other
mathematicians, including Abresch \cite{Abresch}, Bobenko \cite{Bobenko}, Pinkall and Richter
\cite{KamberovGANG}, and Taimanov \cite{Taimanov}.  It is satisfied for a minimal surface ($H = 0$) if and
only if $\psi = (s_1,\,s_2)$ is meromorphic.

\section{Regular homotopy classes and spin structures}\label{topspin}

Let $X_1,X_2:M\to \bbR^3$ be two immersions of a surface into space.
Recall the distinction between regular homotopy equivalence of the
immersions $X_1$, $X_2$, and regular homotopy equivalence of the
corresponding immersed surfaces --- these immersed surfaces are
regularly homotopic if there is a diffeomorphism $h$ of $M$ such
that $X_2$ is regularly homotopic to $X_1\circ h$ --- so this
latter equivalence relation is coarser.

\begin{theorem}\label{homotopy theorem}
Let $X_1,X_2:M\to \bbR^3$ be two immersions of a surface into space, 
let $S_1$, $S_2$ the spin structures induced as in
Theorem~\ref{induced spin structure}, and let $q_1$, $q_2$ be the 
associated quadratic forms as in Theorem~\ref{quadratic theorem}. Then
\begin{enumerate}
\item $X_1$ and $X_2$ are regularly
homotopic if and only if $q_1\equiv q_2$ \mbox{(mod 2)}.
\item The immersed surfaces $X_1(M)$ and $X_2(M)$ are regularly homotopic
if and only if {\em $\Arf q_1 = \Arf q_2$}.
In particular, an immersed surface is regularly homotopic to an
embedding if and only if its $\Arf$ invariant equals $+1$.
\end{enumerate}
\end{theorem}

\begin{pf*}{\it Sketch of proof}
Define $\tilde{q}(\al)$ as half the linking number (mod 2) of the
boundary curves of the image of a tubular neighborhood of $\al$ in
$\bbR^3$. Then (for $i = 1,\,2$)
\[
q_i(\al)=0 \iffspace 
\parbox[c]{1.5in}{\small the Darboux frame along $\al$ is nontrivial 
as an element of $\pi_1(\SO(3))$} \iffspace \tilde{q}_i(\al) = 0.
\]
Hence $q_i\equiv \tilde{q}_i$ (mod 2).
But $X_1$, $X_2$ are regularly homotopic if and only if 
$\tilde{q}_1\equiv\tilde{q}_2$ (mod 2), and the corresponding
immersed surfaces are regularly
homotopic if and only if $\Arf \tilde{q}_1 = \Arf \tilde{q}_2$
(see \cite{Pinkall}).
\end{pf*}

\section{Spin structures and even-order differentials}\label{holo}

Theorem~\ref{holo spin theorem} ties the notion of spin structure with
other concepts from algebraic geometry.
Recall that a {\em theta characteristic} on a Riemann surface
is a divisor $D$ such that $2D$ is the
canonical divisor.

\begin{theorem}\label{holo spin theorem}
Given a Riemann surface $M$,
there are natural bijections between the following sets of objects:
\begin{enumerate}
\item the spin structures on $M$;
\item the complex line bundles $S$ on $M$ satisfying $S\tensor S\cong K$;
\item the theta characteristics on $M$;
\item the classes of non-identically-zero 
meromorphic differential forms on $M$ whose zeros and poles
have even orders, under the equivalence
\[
\eta_1 \sim \eta_2 \iffspace 
  \parbox[c]{2in}{\small $\eta_1/\eta_2 = h^2$
  for some meromorphic function $h$ on $M$.}
\]
\end{enumerate}
\end{theorem}

\begin{pf}
{\it (i)$\iff$(ii).} Given a line bundle $S$ on $M$ satisfying 
$S\tensor S\cong K$, $S$ is a spin structure with squaring map 
$\mu:S\to S\tensor S$ defined by $\mu(s) = s\tensor s$.
Conversely, given a spin structure $S$ on $M$, the bundle map
$\mu(s)\mapsto s\tensor s$ is well-defined and a vector-bundle isomorphism,
so $K$ is isomorphic to $S\tensor S$.

{\it (ii)$\iff$(iii).}
Via the natural correspondence between the line bundles on $M$ with
the divisor classes,
this set of line bundles is bijective with
the theta characteristics.

{\it (iii)$\iff$(iv).}
Again, there is a natural bijection between the
meromorphic differentials with zeros and poles of even orders 
and the theta characteristics. Given such a differential
$\eta$, the corresponding theta characteristic is
$\frac{1}{2}(\eta)$. Moreover,
two such differentials correspond to theta characteristics in
the same linear equivalence class if and only if their
ratio is the square of a meromorphic function on $M$.
For 
$\eta_1/\eta_2 = h^2$ if and only if  
$\textstyle \frac{1}{2}(\eta_1)-\frac{1}{2}(\eta_2) = (h).$
\end{pf}

The spin structures on a compact Riemann surface are also bijective with
the various translates $\vartheta[{a_0\atop b_0}]$ of the theta functions 
on the surface (see \cite{Mumford} for the definition of 
$\vartheta[{a_0\atop b_0}]$).

\section{Spin structures on tori}\label{spin on tori}

We compute the four spin structures on a Riemann torus $T$ together with
their values of $q$.
Let $\bbC/\{2\om_1,2\om_3\} = \hbox{Jac}(T)$ be the Jacobian for $T$, and
let $e_i = \wp(\om_i)$ ($i$ = 1, 2, 3), where $\om_2 = \om_1+\om_3$.
Then $h(u) = (\wp(u),\wp'(u))$ is a conformal diffeomorphism from the Jacobian
to the Riemann surface $M$ defined by $w^2 = 4(z-e_1)(z-e_2)(z-e_3)$. It is
then elementary to show that the four differentials
\vspace{-0.1in}
\[
\begin{array}{rcl}
du &=& dz/w,\\
(\wp(u) - e_i)du &=& (z - e_i)dz/w
\end{array}
\]
define the four distinct spin structures as in Theorem~\ref{holo spin theorem}.
\noindent
With $\al_i$ the generator of $H_1(T,\bbZ_2)$ defined by
$\al_i:[0,1]\to\hbox{Jac($T$)}$, $\al_i(t) = 2 t \om_i$,
the values of $q$ and $\Arf q$ are tabulated.

\vspace{-0.1in}
\begin{table}[htbp]
\centering
\caption{Values of $q$ and $\Arf q$ for spin structures on tori}
\vspace{-.12in}
$
\begin{array}{c||c|c|c|c||c}
          \eta & 
      q_\eta(0) & q_\eta(\al_1) & q_\eta(\al_2) & q_\eta(\al_3)
      & \Arf q_\eta\\ \hline
            du & 0 & 1 & 1 & 1 & -1\\
(\wp(u)-e_1)du & 0 & 1 & 0 & 0 & +1\\
(\wp(u)-e_2)du & 0 & 0 & 1 & 0 & +1\\
(\wp(u)-e_3)du & 0 & 0 & 0 & 1 & +1
\end{array}
$
\end{table}
An immersion corresponding to $q$ for which $\Arf q = +1$ 
is regularly homotopic to the torus standardly embedded in $\bbR^3$.
The value
$\Arf q = -1$ corresponds to the twisted torus, which can be realized as
the ``diagonal'' double covering of the standardly embedded torus as shown,
but is not regularly homotopic to an embedding.

The more general case of hyperelliptic Riemann surfaces is considered in
Appendix \ref{hyper appendix}, where the spin structures and their
corresponding quadratic forms are computed explicitly.
\vspace{-.1 in}
\begin{figure}[htbp]
\setlength{\unitlength}{0.6cm}
\centering\begin{picture}(4,3.3)(0,0.8)
\small
\multiput(1,1)(0,2){2}{\line(1,0){2}}
\multiput(1,1)(2,0){2}{\line(0,1){2}}
\multiput(2,0)(2,2){2}{\line(-1,1){2}}
\multiput(2,0)(-2,2){2}{\line(1,1){2}}

\linethickness{.005cm}
\multiput(1,1)(0,0.1){20}{\line(1,0){2}}
\multiput(1,1)(0.1,0){20}{\line(0,1){2}}

\put(4,3){\vector(-2,-1){1.7}}
\put(3.7,3.2){standard torus}

\end{picture}
\vspace{.13in}
\caption{The twisted torus}
\end{figure}
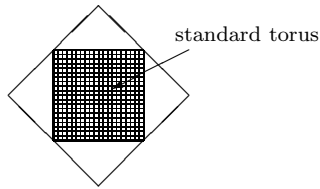

\section{Group action on spinors}\label{group action}

The automorphism group of $Q$ is the linear conformal group
\[
\bbC^*\times\SO(3,\bbC).
\]
The orbit of a conformal immersion $X$ with $\partial X = \om \in Q$ under this action is an 8-real-dimensional 
family of immersions. (This action, however, will not respect the
vanishing of periods --- see Section \ref{periods}.) The subgroup 
\[
\bbR^+\times\SO(3,\bbR)
\]
is the group of similarity transformations of Euclidean 3-space.
Hence the homogeneous space
\begin{equation}\label{homog1}
\left(\bbC^*\times\SO(3,\bbC)\right)/
       (\bbR^+\times\SO(3,\bbR)) \cong
  S^1\times\left(\SO(3,\bbC)/\SO(3,\bbR)\right).
\end{equation}
is the 4-real-dimensional parameter space of non-similar surfaces in
the above orbit. 

In terms of spinors, we use the two-fold spin covering group $\GL(2,\bbC)$ to get this family of immersions,
as justified by the following well-known fact
(see, for example, \cite{Gilbert}, \cite{Rees}). Some details are given in
Appendix \ref{appendix group action}.

\begin{theorem} \label{two-fold cover}
There is a unique two-fold covering homomorphism (spin)
{\em
\[
T:\GL(2,\bbC)\to \bbC^*\times \SO(3,\bbC)
\]
}
such that for any {\em $A\in\GL(2,\bbC)$},
\begin{equation}\label{T commutes}
T(A) \sig = \sig A,
\end{equation}
where $\sig:\bbC^2\to Q$ is as in equation (\ref{sigma}), and
$A$ and $T(A)$ act by left multiplication on $\bbC^2$ and $Q$ 
respectively.
Moreover, $T$ restricts to a two-fold covering 
of {\em $\SL(2,\bbC)$}, {\em $\bbR^*\times\SU(2)$}, and
{\em $\SU(2)$} onto {\em $\SO(3,\bbC)$}, {\em $\bbR^*\times\SO(3)$}, and
{\em $\SO(3)$}.
\end{theorem}

Lifting the group action on $Q$ to an action on $\bbC^2\setminus\{0\}$ 
via $T$,
the homogeneous space in equation (\ref{homog1}) can also be written
\begin{equation}\label{homog2}
  \left(\GL(2,\bbC)\right)/
  \left(\bbR^*\times\SU(2)\right)\cong
  S^1\times\left(\SL(2,\bbC)/\SU(2)\right)\cong
  S^1\times{H^3},
\end{equation}
where $H^3$ is hyperbolic three-space.
The $S^1$ factor gives rise to the well-known ``associate family'' of
minimal surfaces, which are locally isometric and share a common Gauss
map. The other factor has a simple (though apparently less known)
geometric interpretation as well. The Gauss map is the ratio of two
spinors, so $\SO(3,\bbC) \cong \PSL(2,\bbC)$ acts on the Gauss map via
post-composition with a fractional linear transformation of $S^2$;
indeed, the quotient by $\SO(3,\bbR) \cong \PSU(2)$ leaves $H^3$,
so the second factor can be thought of
as the non-rigid M\"obius deformations of the Gauss map.

\section{Periods}\label{periods}

Given an immersion $X:M\to \bbR^3$, the {\em period} around a 
simple closed curve $\gam\subset M$ is the vector in $\bbC^3$
\[
\int_{\gam}\partial X.
\]
If the real part of a period is not $(0,0,0)$, the resulting surface is
periodic and does not have finite total curvature. It is a
considerable problem
to ``kill the periods'' --- that is, choose parameters so that the
integrals around every simple closed curve in $M$ generates purely imaginary period vectors.
Non-zero periods can arise along two kinds of simple closed curves:
\begin{enumerate}
\item[$\bullet$] a simple closed curve around an end $p\in M$,
\item[$\bullet$] a non-trivial simple closed curve in $H_1(M,\bbZ).$
\end{enumerate}
For minimal surfaces, $\partial X$ is meromorphic, and a simple closed
curve $\gam$ around an end $p\in M$ has period
\[
\int_{\gam} \partial X = 2\pi i \res_p\partial X.
\]
This integral is zero for minimal surfaces with embedded planar ends (see Section \ref{ends}).

Using the spinor representation,
the condition that a period along a closed curve $\gam\subset M$
be pure imaginary can be expressed by
\[
\int_\gam \left(s_1^2-s_2^2,i(s_1^2+s_2^2), 2s_1 s_2\right)\in i\bbR^3,
\]
equivalent to
\begin{equation}\label{period1}
\disp\int_{\gam} s_1^2 = \disp\ov{ \int_{\gam}s_2^2}\andspace
\disp\int_{\gam} s_1 s_2 \in i\bbR.
\end{equation}
These equations are preserved by the group $\bbR^*\times\SU(2)$
of homotheties.

\section{Spin structures and nonorientable surfaces}\label{nonorientable}

To deal with immersions of a nonorientable surface $M$ into space, we
pass to the oriented two-fold cover of $M$. The following rather
technical results are required in Part III. Without proof we state:

\begin{lemma}\label{antipodal}
Let
{\em
\[
\begin{array}{l}
A:S^2\to S^2 \hbox{\ {\it be the antipodal map,}}\\
A_*:T(S^2)\to T(S^2) \hbox{\ {\it the derivative of $A$,}}\\
\hat{A}_*:\Spin(S^2)\to\Spin(S^2) 
    \hbox{\ {\it one of the lifts of $A_*$ to $\Spin(S^2)$ .}}
\end{array}
\]
}
Then, in the coordinates of Section \ref{spinor representation}, we have
{\em
\[
A_* = \conj\hspace{0.5in}
\hat{A}_* =
\pm
\left(\begin{array}{cc} 0 & i\\ -i & 0\end{array}\right)
\circ\conj.
\]
}
\end{lemma}

\begin{figure}[htbp]\label{figure:antipodal}
\setlength{\unitlength}{1.8cm}
\centering\begin{picture}(2,2.5)(-.3,0)
\small
\multiput(0.1,0)(0,1){3}{\vector(1,0){0.7}}
\multiput(-.3,0.8)(1.4,0){2}{\vector(0,-1){0.6}}
\multiput(-.3,1.8)(1.4,0){2}{\vector(0,-1){0.6}}

\put(-0.4,-0.08){$S^2$}
\put(1,-0.08){$S^2$}
\put(-0.6,.9){$T(S^2)$}
\put(1,.9){$T(S^2)$}
\put(-.8,2){$\Spin(S^2)$}
\put(1,2){$\Spin(S^2)$}

\put(.35,0.1){$A$}
\put(.35,1.1){$A_*$}
\put(.35,2.1){$\hat{A}_*$}
\put(1.2,1.5){$\sig$}

\end{picture}
\caption{Lifts of the antipodal map}
\end{figure}


\begin{theorem}\label{nonorientable theorem}
Let $M$ be a nonorientable Riemann surface, and
$X:M\to\bbR^3$ a conformal immersion of $M$ into space.
Let $\pi:\widetilde{M}\to M$ be an oriented
double cover of $M$, and $\widetilde{X}=X\circ\pi$ the lift of
$X$ to this cover. Let $I:\widetilde{M}\to \widetilde{M}$ the order-two
deck transformation for the cover. With $\om = \partial\widetilde{X}$, 
and in the notation of Lemma~\ref{antipodal},
we have
\begin{enumerate}
\item $gI = Ag$,
\item $\om I_* = A_*\om$,
\item $\psi\hat{I}_* = \pm\hat{A}_*\psi$.
\end{enumerate}
\end{theorem}

%

We remark that the proper treatment of nonorientable surfaces should
really be via ``pin'' structures (Pin($n$) being the corresponding
two-sheeted covering group of O($n$)), and that in this case we should
have an analytic formula (in analogy with that in Appendix B for the $\bbZ_2$-valued Arf invariant on hyperelliptic surfaces) for the full $\bbZ_8$-valued Arf invariant
\[
\frac{1}{\sqrt{\# H}}\sum_{ \al\in H} i^{q(\al)}
\]
of the associated $\bbZ_4$-valued quadratic form $q$ on
$H =H^1(M,\bbZ_2)$.
\nopagebreak
%
%

\typeout{_______________________________________________ part2.tex}
\part{Minimal Immersions with Embedded Planar Ends}

\label{chapter:ends}
The first section of this part of our paper discusses the behavior of
a minimal immersion at an embedded planar end.  Lemma~\ref{lemma:ends}
translates this geometric behavior to a necessary and sufficient
algebraic condition on the order and residue of the immersion at the
end. Arising naturally from this algebraic condition is a certain
vector subspace $\calK$ of holomorphic spinors which generates all
minimal surfaces with embedded planar ends (Theorem~\ref{theorem: moduli
space}).  More precisely, two sections chosen from $\calK$ form
the spinor representation of a minimal surface, and conversely, any
such surface must arise this way (Theorem~\ref{K condition}).
However, such a surface is usually periodic, and possibly a branched
immersion.  In order to compute $\calK$ explicitly, a skew-symmetric
bilinear form $\skewform$ is constructed (Definition
\ref{skewform}) whose kernel is closely related to the space $\calK$.
In Part III, this form is used to prove existence and non-existence
theorems for a variety of examples.

\section{Algebraic characterization of embedded planar ends}\label{ends}

The geometric condition that an end of a minimal immersion be embedded
and planar can be translated to algebraic conditions (see, for example,
\cite{Callahan}).  Let $X:D\setminus\{p\}\to\bbR^3$ be a conformal
minimal immersion of an open disk $D$ punctured at $p$ such that
$\lim_{q\rightarrow p}|X(q)|=\infty$.  The image under $X$ of a small
neighborhood of $p$ (and by association, $p$ itself) is what we shall
refer to as an {\em end}.  The behavior of the end is determined by
the residues and the orders of the poles of $\partial X$ at $p$ as follows.

Let $\zeta_1$, $\zeta_2$, $\zeta_3$ be defined by
\[
\partial X = (\zeta_1 - \zeta_2, i(\zeta_1+\zeta_2), 2\zeta_3).
\]
The Gauss map for this immersion (see \cite{Osserman}) is
\[
 g = \zeta_3/\zeta_1 = \zeta_2/\zeta_3.
\]

First note that for $X$ to be well-defined, we must have for any closed 
curve $\gam$, which winds once around $p$,
\[
0=\re \int_\gam \partial X = \re (2 \pi i \res_p \partial X),
\]
and so $\res_p\partial X$ must be {\em real}.  Assume this, and
assume initially that the limiting normal to the end is upward (that
is $g(p)=\infty$).  In this case,
\[
\ord_p\zeta_2 < \ord_p\zeta_3 < \ord_p\zeta_1,
\]
so the first two coordinates of $X(q)$ grow faster than does the third
as $q \rightarrow p$.

It follows that $\ord_p\zeta_2$ cannot be $-1$, because if it were then
\[
\res_p\partial X = (-\res_p\zeta_2,i\res_p\zeta_2,0)
\]
would not be real.  Hence $\ord_p\zeta_2 \le -2$.  The image under
$X$ of a small closed curve around $p$ is a large curve which winds around the
end $|\ord_p\zeta_2|-1$ times.  The end is embedded precisely when 
$\ord_p\zeta_2=-2$. 

If an end is embedded, its behavior is determined by the vanishing or
non-vanishing of the residues of $\partial X$.  For an embedded end,
$-2 = \ord_p\zeta_2 < \ord_p\zeta_3$, so $\zeta_3$ has either a simple
pole or no pole.  If $\zeta_3$ has a simple pole (and hence also a
residue), the end grows logarithmically relative to its horizontal
radius and is a {\em catenoid} end. If $\zeta_3$ has no pole, the end is
asymptotic to a horizontal plane and is called a {\em planar} end.
Moreover, in this latter case, $\res_p\zeta_2$ must vanish (again, if
it did not, $\res_p\partial X$ would not be real), and so
$\res_p\partial X = (0,0,0)$.

For an end in general position the same conclusions hold,
because a real rotation affects neither $\ord_p\partial X$ nor
the reality or vanishing of $\res_p \partial X$.  In summary, 
we have

\begin{lemma}\label{lemma:ends}
Let $X:D\setminus\{p\}\to\bbR^3$ be a conformal minimal immersion of
a punctured disk.
Then $p$ is an embedded planar end if and only if
{\em
\[
\ord_p \partial X = -2 \andspace \res_p\partial X = 0,
\]
}
where {\em $\ord_p\partial X$} denotes the minimum order at $p$ of the
three coordinates of $\partial X$.
\end{lemma}

\section{Embedded planar ends and spinors}\label{endspin}

The conditions in the lemma above can be translated into conditions on
the spinor representation of the minimal immersion. This leads to the
definition of a space $\calK$ of spinors, pairs of which form the spinor
representation satisfying the required conditions.

Throughout the rest of Part II, the following notation is fixed:
\begin{equation}\label{setup}
\begin{array}{l}
\hbox{$M$ is a compact Riemann surface of genus $g$,}\\
\hbox{$K = T(M)$ is the canonical line bundle,}\\
\hbox{$S$ is a spin structure on $M$,}\\
\hbox{$P = [p_1] + \dots + [p_n]$ is a divisor of $n$ distinct points.}
\end{array}
\end{equation}
The points $\upto{p}{1}{n}$ will eventually be the ends of a minimal
immersion of $M$ whose spinor representation will be a pair of
sections of $S$.

Let $H^0(M,\calO(S))$ and $H^0(M,\calM(S))$
denote respectively the vector spaces of holomorphic and 
meromorphic sections of $S$.
Define
\begin{equation}\label{spaces}
\begin{array}{lclcl}
\calF \!\!\!\andeq\!\!\! \calF_{M,S,P}\!\!\!\andeq\!\!\!\{s\in H^0(M,
\calM(S))\suchthat (s)\ge -P\}\\
\calH \!\!\!\andeq\!\!\!\calH_{M,S}\!\!\!\andeq\!\!\! H^0(M,\calO(S))\\
\calK \!\!\!\andeq\!\!\!\calK_{M,S,P}\!\!\!\andeq\!\!\! 
  \{s\in\calF\suchthat \ord_ps\ne 0 \hbox{ and }
  \res_ps^2 = 0 \hbox{ for all $p\in\supp P$}\}.
\end{array}
\end{equation}

We remark that
\begin{equation}\label{calK}
s\in\calK \iffspace
   \parbox[c]{\boxwidth}{\small the constant term in the expansion of
                                $s$ vanishes at each $p\in P$.}
\end{equation}
Thus $\calH$ and $\calK$ are linear subspaces of $\calF$. 

\begin{theorem}\label{K condition}
Let $X:M\to\bbR^3$ be a minimal immersion 
with spinor representation $(s_1,s_2)$. Then
$p\in M$ is an embedded planar end if and only if
$s_1$, $s_2\in\calK$ and
at least one of $s_1$, $s_2$ has a pole at $p$.
\end{theorem}

\begin{pf}
By Lemma~\ref{lemma:ends}, $p$ is an embedded planar end if and only if
\[
\ord_p\partial X = -2 \andspace \res_p\partial X = 0.
\]
The first of these equations is equivalent to the condition 
\[
\hbox{$s_1$, $s_2\in\calF$, and at least one of $s_1$, $s_2$ has a 
     pole at $p$.}
\]
Suppose that this condition is satisfied, and let $z$ be a conformal
coordinate near $p$ with $z(p) = 0$, $\varphi^2 = dz$ and
\[
s_1 = \left(\frac{a_{-1}}{z} + a_0 + \dots\right)\varphi
\andspace
s_2 = \left(\frac{b_{-1}}{z} + b_0 + \dots\right)\varphi
\]
be expansions of $s_1$ and $s_2$ respectively.
Then the condition $\res_p\partial X = 0$ is equivalent to
\begin{eqnarray*}
\res_p s_1^2   \andeq 2 a_{-1} a_0 = 0 \\
\res_p s_1 s_2 \andeq a_{-1}b_0 + a_0 b_{-1} = 0\\
\res_p s_2^2    \andeq 2 b_{-1} b_0 = 0
\end{eqnarray*}
or, under the assumption that $s_1$ or $s_2$ has a pole at $p$ (i.e.,
$a_{-1} \ne 0$ or $b_{-1} \ne 0$)
\[
a_0 = b_0 = 0.
\]
This is to say that $s_1$, $s_2 \in \calK$.
\end{pf}

\section{Moduli spaces of minimal surfaces with embedded planar ends}
The spinor representation yields the following general result for a fixed spin structure $S$ over a fixed Riemann surface $M$
with fixed punctures at $P$.

\begin{theorem}\label{theorem: moduli space}
Let $P$ be a divisor and $S$ a spin structure on $M$
as in equation (\ref{setup}), and let
$\calK = \calK_{M,S,P}$ as in equation (\ref{spaces}) with
$m = \dim\calK\ge 2$. Then the space of complete minimal (branched,
possibly periodic) immersions of $M$ into $\bbR^3$ with finite total
curvature and embedded planar ends at {\em $\supp P$} is the complex
$2(m-1)$-dimensional manifold {\em $\Gr_{\!2}(\calK)\times\left(S^1\times H^3\right)$}.  All these immersions are regularly homotopic and in the class determined
by $S$.
\end{theorem}
\par\begin{pf}
Fix a point of the Grassmanian $\Gr_{\!2}(\calK)$, which
represents a two-dimensional complex plane in $\calK$.  Let $(t_1, t_2)$ be a basis for
this plane.  Then
the family of (branched, possibly periodic) minimal immersions is given by $X=\re F$, where
\[
F = \int(s_1^2-s_2^2, i(s_1^2+s_2^2), 2 s_1 s_2)
\]
and 
\[
\left(\begin{array}{c}s_1\\s_2\end{array}\right)
=
R
\left(\begin{array}{c}t_1\\t_2\end{array}\right)
\]
for some $R\in\GL(2,\bbC) = \bbC^*\times\SL(2,\bbC)$.
The surfaces are identical (up to a rotation or dilation in space)
when $R\in\bbR^*\times\SU(2)$.  Thus a parameter space for this family
of surfaces is $S^1\times H^3$ (see Section \ref{group action}).
\end{pf}

We remark that $\Gr_{\!2}(\bbC^m)\times\left(S^1\times H^3\right)$ is
actually a quaternionic manifold.

\section{The vector spaces $\calF$, $\calH$ and $\calK$}\label{F dimension}

When the number of ends exceeds the genus, the dimension of the space
$\calF$ of meromorphic sections of $S$ (with at most simple poles at
the ends) is computable, and any holomorphic section in $\calK$ must vanish.

\begin{theorem}\label{p:space}
Let $M$, $P$, and $S$ be as in equation (\ref{setup}),
and $\calF$, $\calH$, and $\calK$ as in equation (\ref{spaces}).
Let {\em $g=\genus(M)$} and {\em $n = \deg P$}. Then, under the
assumption that $n \ge g$,
\begin{enumerate}
\item $\dim\calF = n$;
\item $\calK \cap \calH=0$.
\end{enumerate}
\end{theorem}

\begin{pf*}
{\it Proof of (i)}
The dimension of $\calF$ can be computed by means of the 
Riemann-Roch theorem (see, for example, \cite{Gunning1}) which states
\[
\dim H^0(M,L) - \dim H^0(M,K\tensor L^*) = \deg L - g + 1
\]
for an arbitrary line bundle $L$.
Let $R$ be the line bundle corresponding to the divisor $P$, and let
$L = S\tensor R$.  Then:
\begin{enumerate}
\item[$\bullet$] 
$H^0(M,L) \cong \calF$ by the isomorphism  $s\tensor r \mapsto s$, 
where $r$ is a section of $R$ with divisor $P$;
\item[$\bullet$] 
$H^0(M,K\tensor L^*)=0$, since $\deg(K\tensor L^*) = g-1-n$,
which is negative by hypothesis;
\item[$\bullet$] 
$\deg L - g + 1 = n$;
\end{enumerate}
from which it follows that
\[
\dim\calF = \dim H^0(M,L) = n.
\]

{\it Proof of (ii).}
Let $s\in\calK\cap\calH$ be a section which is not identically zero.  
Since $s\in\calK$, we have that $\ord_p s\ne 0$ for all $p\in\supp P$ ---
that is, at each such $p$, $s$ has either a pole or a zero.  But since 
$s\in\calH$, $s$ cannot have a pole at $p$, and hence has a zero, so
$(s)\ge P$.  Conversely, if $(s)\ge P$, then $s\in\calK\cap\calH$, so
\[
\calK\cap\calH = \{s\in\calF\suchthat (s)\ge P\}.
\]
Thus for $s\in \calK\cap\calH$ not identically zero,
\[
n \le \deg s = g-1.
\]
Hence if $n\ge g$, then $\calK\cap\calH=0$.
\end{pf*}

\section{A bilinear form $\skewform$ which annihilates $\calK$}\label{bilinear}

To understand the vector space $\calK=\calK_{M,S,P}$ more explicitly,
a skew-symmetric bilinear form $\skewform$ is defined on $\calF$ whose
kernel contains $\calK$.  This form may then be used in many cases to
compute $\calK$, and thereby moduli spaces of minimal surfaces with
embedded planar ends.  Since, in effect, we are now interested in
varying $M$ and $P$, it is natural to consider quadratic differentials
on $M$ with poles at $P.$  If $\Phi$ is such a differential (that is, a
meromorphic section of $K\tensor K$) with expansion
\[
 \Phi = \sum_{k = -\infty}^\infty a_k z^k dz^2 
\]
at a point $p$ on $M$ in the conformal
coordinate $z$ with $z(p) = 0,$ then the number $a_{-2}$ is independent of
this choice of coordinate, and is denoted in what follows by
$\qres_p \Phi.$  We shall use the Hopf differential $\Phi =
s_1 \partial s_2 - s_2 \partial s_1$ from Section \ref{spinor representation}.

\begin{definition}\label{skewform}
With $M$, $P = \sum_{k=1}^n [p_k]$ and $S$ as in equation
(\ref{setup}) define
$\skewform = \skewform_{M,P,S}:\calF \times \calF \to \bbC$ by
{\em
\[ 
\skewform(s_1,s_2) =
-\frac{1}{2}\sum_{p \in P} \qres_p (s_1 \partial s_2 - s_2 \partial s_1),
\]
}
and define $\widehat{\skewform} : \calF \to \calF ^{*}$
by $(\widehat{\skewform}(s))(t) = \skewform(s, t).$
\end{definition}

\begin{theorem}\label{p:skew}
With $\calH$, $\calK$ as in equation (\ref{spaces}),
$\skewform$ satisfies the following:
\begin{enumerate}
\item $\skewform$ is a skew-symmetric bilinear form on $\calF$;
\item  $\ker \widehat{\skewform}\supseteq\calK + \calH$;
\item if $\calK \cap\calH=0$, then $\ker \widehat{\skewform} = 
         \calK \oplus \calH$;
\item if {\em $n=\deg P\ge g=\genus (M)$}, then
  $\ker\widehat{\skewform} = \calK \oplus \calH$.
\end{enumerate}
\end{theorem}

\begin{pf}
Part (i) is immediate from the definition of $\skewform.$

{\it (ii) and (iii).}  We use a choice of coordinates to factor
$\widehat{\skewform}$
into a pair of linear maps whose kernels are $\calH$ and $\calK$ respectively.

For each $k \in \{1, \dots, n\},$ let $z_k$ be a conformal coordinate
in a neighborhood $U_k$ of $p_k$ with $z_k(p_k) = 0.$ Let $\varphi_k$ be
a spinor on $U_k$ with $\varphi_k^2 = dz_k.$ With these choices, for any
spinor $s,$ let $\al_r^k(s)$ denote the coefficient of $z_k^r$ in the
local expansion of $s/\varphi_k$ at $p_k$.  That is, the expansion of $s$
at $p_k$ is
\[
s =   \left(
  \frac{\al_{-1}^k(s)}{z_k} + \al_0^k(s) +
      \dots
  \right)\varphi_k.
\]
Then
\[
0 = \sum_{p \in P} \res_p s_1 s_2 = \sum_{k=1}^{n}
  \left(\al_{-1}^k(s_1) \al_0^k(s_2) +
  \al_0^k(s_1) \al_{-1}^k(s_2)\right)
\]
and so
\begin{eqnarray}\label{local}
\skewform(s_1, s_2) \andeq
  -\frac{1}{2} \sum_{k=1}^n (\al_{-1}^k(s_1) \al_0^k(s_2) -
   \al_0^k(s_1) \al_{-1}^k(s_2)) \\
\andeq \sum_{k=1}^n \al_0^k(s_1) \al_{-1}^k(s_2) =
                    - \sum_{k=1}^n \al_{-1}^k(s_1) \al_0^k(s_2).
\end{eqnarray}

Let the linear maps
$A_r:\calF \to \bbC^n$ ($r = -1, \ 0$) be defined by
\[
A_r(s) = (\al_r^1(s), \dots, \al_r^n(s)).
\]
Then, identifying $(\bbC^n)^*$ with $\bbC^n$ in the natural way,
$\skewform$ factors as
\[
\widehat{\skewform} = (A_{-1})^* \circ A_0 = -(A_0)^* \circ A_{-1}.
\]
Part (ii) then follows from the facts that
\[
\calH = \ker A_{-1} \andspace \calK = \ker A_0,
\]
the latter by equation (\ref{calK}).

Finally, since for any linear maps
$X \stackrel{A}{\to} Y \stackrel{B}{\to} Z,$
\[
\dim(\ker B \circ A) =
  \dim(\ker A) + \dim(\image A \cap \ker B) \le \dim(\ker A) + \dim(\ker B),
\]
one has
\[
\dim(\ker \widehat{\skewform}) \le \dim(\ker A_{-1}) + \dim(\ker A_0) =
  \dim(\cal H) + \dim(\calK).
\]
It follows, under the assumption that $\calH\cap\calK = 0,$ that
$\ker \widehat{\skewform} = \calH \oplus \calK.$  This proves part (iii).

{\it (iv).}  This follows directly from part (iii)
and Theorem~\ref{p:space}(ii).
\end{pf}

%
%

\typeout{_______________________________________________ part3.tex}
\part{Classification and Examples}

In the first half of Part III, the skew-symmetric form $\skewform$
developed in Part II is used to investigate minimal genus zero
surfaces with embedded planar ends.  The first two sections
demonstrate the non-existence of examples with 2, 3, 5, or 7
ends, and the dimension of the moduli space of examples with 4, 6,
8, 10, 12 and 14 ends is computed. The following two
sections compute explicitly the moduli spaces for the families with
4 and 6 ends, and in section \ref{projective}, the moduli space
of the three-ended projective planes is investigated.

The remaining sections are devoted to minimal immersions in the
regular homotopy classes of tori and Klein bottles with embedded
planar ends. In Sections \ref{du} and \ref{nondu}, the skew-symmetric
form $\skewform$ is computed for the twisted and the untwisted tori.
This computation is then used to show the nonexistence and existence
of various examples. In Section \ref{3ends} it is shown that no such
tori exist with three ends, and in Section \ref{4ends2}, is found a
real two-dimensional family of twisted immersions with four ends
exists on each conformal type of torus. An amphichiral minimal Klein bottle with four
embedded planar ends is constructed in Section \ref{kleinex}.

All of these surfaces are found (or shown not to exist) by the
following general method: after computing $\skewform$ on a simple basis,
its pfaffian, which is a function of the ends, is set to zero. The
resulting condition on the placement of the ends --- that is, the
determinantal variety --- together with further conditions arising from
the demand that the immersion have no periods and no branch points,
forms a set of equations whose simultaneous solution (or impossibility
of solution) gives the desired result.

\section{Existence and non-existence of genus-zero surfaces}\label{nonexist}

The non-existence of genus zero minimal unbranched immersions
with 3, 5 or 7 embedded planar ends was first proved in a case-by-case
manner in \cite{Bryant2}. The following is a new proof, using the ideas
of Part \ref{chapter:ends}.

\begin{theorem}\label{theorem:nogenus0}
There are no complete minimal branched or unbranched immersions of a
punctured sphere into space with finite total curvature and $2$, $3$,
$5$, or $7$ embedded planar ends. There exist unbranched examples with
$4$, $6$, and any $n\ge8$ ends.
\end{theorem}

\begin{pf}
Examples with $2p$ ends ($p\ge 2$) are given in \cite{Kusner2}, and
with $2p+1$ ends ($p\ge 4$) in \cite{Peng}. For the cases $n=3$, $5$,
or $7$, by Lemma~\ref{dimension lemma}(ii) below,
$2\le \dim\calK\le [\sqrt{n}] \le 2$ (here $[q]$ denotes the greatest
integer less than or equal to $q$), so $\dim\calK=2$, which contradicts
Lemma~\ref{dimension lemma}(i) that $n-\dim \calK$ is even. The case
$n=2$ is proved in \cite{Kusner2} (or is proved likewise by
Lemma~\ref{dimension lemma} below).
\end{pf}

We remark that there is also a simple topological proof of the
non-existence of genus zero examples with 3 ends, using ideas in
\cite{Kusner1} and \cite{Kusner3}.  The trick is to exploit the
$\SO(3,\bbC)$-action discussed in Section \ref{group action}
to deform the Gauss map --- on a punctured sphere with planar ends
there is no period obstruction to doing this --- so the
compactified $S^2$ is in general position with a unique (tranverse) triple-point,
which is impossible.  (By carefully treating the periods introduced by
this explicit $\SO(3,\bbC)$ deformation of the Gauss map, the same kind
of argument should generalize to exclude orientable minimal surfaces of
any genus with three embedded planar ends --- see Section \ref{3ends}
for a different proof in case of tori.)

\begin{lemma}\label{dimension lemma}
Let $P$ be a divisor on the Riemann sphere $S^2$ as in 
equation (\ref{setup}) with $n=\deg P\ge 2$, and
let $\calK = \calK_{S^2,S,P}$ be as in equation (\ref{spaces}).  Then
\begin{enumerate}
\item $n-\dim\calK$ is even;
\item If there exists a complete branched or unbranched
minimal immersion of $S^2$ into space with finite total
curvature and $n$ embedded planar ends in {\em $\supp (P)$}, then 
$2 \le \dim\calK \le \sqrt{n}.$
\end{enumerate}
\end{lemma}

\begin{pf*}
{\it Proof of (i)}
By Theorem~\ref{p:skew}, $\ker\skewform = \calK\oplus\calH$.
But $\calH=0$ because there are no holomorphic differentials on the sphere,
so $\ker\skewform = \calK$.  Since $\skewform$ is skew-symmetric,
$\rank\skewform = n-\dim\calK$ is even (see Appendix \ref{pfaffian}).

{\it (ii).}
The sections $s_1$ and $s_2$ in the spinor representation $(s_1,s_2)$ 
of such a surface are independent, showing the inequality $2\le\dim\calK$.
To show the other inequality,  let $z$ be the standard conformal coordinate
on $S^2=\bbC\cup\{\infty\}$, and let $P=\sum[a_i]$ (where the $a_i\in\bbC$ are distinct)
be the divisor of the
$n$ ends.  Let
$g_1\eta,\dots,g_m\eta$ be a basis for $\calK$, where $\eta^2=dz$.
Define $f:S^2=\CPone\to \bbC\bbP^{m-1}$ by 
\[
f=(\upto{g}{1}{m}).
\]
Then $f$ is well-defined and holomorphic even at the common zeros and the
common
poles of $\upto{g}{1}{m}$.
Let
\[
\textstyle h(z) = \prod(z-a_i).
\]
It follows from
\[
(hg_i) = (h) + (\eta) + (g_i\eta) \ge (P-n[\infty])+[\infty]-P = 
   -(n-1)[\infty]
\]
that
\[
d_0 = \deg f \le n-1.
\]

To show that $f$ ramifies at each $a\in\supp P$,
let $h_i(z) = (z-a)g_i(z)$.
Then $h_i$ does not have a pole at $a$.  Moreover, since by hypothesis
there exists a
minimal surface with ends at $\supp P$, at least one of the $g_i$ has a pole
at $a$, so the $h_i$ cannot all be zero at $a$.
Hence the appropriate condition that  $f$ ramify at $a$ is
\[
\left.\left(h_i h_j' - h_i'h_j\right)\right|_a=0 \hbox{ for all $i$, $j$.}
\]
This is satisfied because of the condition (\ref{calK}) defining $\calK$:
the expansion of $g_i$ at $a$ is
\[g_i = \frac{c_i}{z-a} + \ooo(z-a),
\]
so the expansion of $h_i$ at $a$ is
\[
h_i = c_i + \ooo(z-a)^2,
\]
and so $h_i'(a)=0$ for all i.
Since $f$ ramifies at each $a\in\supp P$, we have
\[
r_0 = \hbox{ ramification index of $f$ }\ge n.
\]

Now let $f_k: \bbC{\bbP}^1\longrightarrow {\bbP}(\Lambda^{k+1}\bbC^m)$ defined by
$f_k = f \wedge f' \wedge \dots \wedge f^{(k)}$ in $\bbC^m$
be the $k^{\makebox{\rm th}}$ associated curve of $f$, and use the 
Pl\"{u}cker formulas (an extension of the Riemann-Hurwitz formula --- see 
\cite{G&H}) which on $\bbC{\bbP}^1$ are
\[
-d_{k-1} + 2 d_k - d_{k+1} - 2 = r_k,
\]
where $d_k$ is the degree of $f_k$, and $r_k$ is the ramification index
of $f_k$.
In the table below,
multiplying the numbers on the left by the inequalities on the
right and adding yields
\[
d_0 \ge (m+n)(m-1)/m.
\]
But $n-1\ge d_0$, so it follows that $n\ge m^2$.
\end{pf*}
\begin{table}[htbp]
\centering
\caption{Values for the Pl\"{u}cker formulas}
\vspace{.15in}
$
\setlength{\arraycolsep}{.3ex}
\begin{array}{c|ccccccccccc}
m-1\;&          & &2d_0    &-&d_1    &-&2&=&r_0    &\ge&n\\
m-2  &      -d_0&+&2d_1    &-&d_2    &-&2&=&r_1    &\ge&0\\
\vdots &&&&\vdots&&&&&\vdots&&\vdots \\
2    &\;-d_{m-4}&+&2d_{m-3}&-&d_{m-2}&-&2&=&r_{m-3}&\ge&0\\
1    &\;-d_{m-3}&+&2d_{m-2}& &       &-&2&=&r_{m-2}&\ge&0\\
\end{array}
$
\end{table}

We may now compute the moduli spaces of genus-zero examples
with an even number of punctures (ends).

\begin{theorem}\label{theorem:small ends}
For each $p\ge 2$ there exists a real $4(p-1)$-dimensional family of
minimal branched immersions of spheres punctured at $2p$ points with finite
total curvature and embedded planar ends.  For $2\le p\le 7$, the moduli
space of such immersions is exactly $4(p-1)$-dimensional.
\end{theorem}

\begin{pf}
Let $P=\sum[a_i]$ be a divisor of degree $2p$ on $S^2$, and $S$ the
unique spin structure on $S^2$.  Let $\calH$ and $\calK$ be as in
equation
(14.11).  Then $\pfaffian\skewform = 0$ (see Appendix \ref{pfaffian})
if and only if 
$\dim\calK\ge 2$ if and only if
there exists a surface with $2p$ ends at $\supp P$.
Counting real dimensions, the space of $2p$ ends is $4p$-dimensional;
the \Moebius transformations of $S^2$
reduce the dimension by $6$, and the pfaffian
condition on the ends reduce the dimension by another $2$, so the 
space of ends which admit surfaces is $(4p-8)$-dimensional.
For each admissible choice of ends, by Theorem \ref{theorem: moduli space}
there is a real $(4\dim\calK -4)$-dimensional space of surfaces.
Altogether, this totals $4p + 4\dim\calK - 12$, which is at least $4p-4$
since $\dim\calK\ge 2$.

In the case that $2\le p\le 7$, by Lemma~\ref{dimension lemma}(ii),
$2\le \dim\calK \le [\sqrt{2p}] \le [\sqrt{14}] = 3$,
so $\dim\calK$, being even, must be exactly 2.
\end{pf}

\section{$\skewform$ on the Riemann sphere}\label{compute}

For the examples in Sections \ref{4ends}--\ref{projective} we need to compute
the skew-symmetric form $\skewform$ from Section \ref{bilinear} on the Riemann sphere.
Let $z$ be the standard conformal coordinate on $S^2=\bbC\cup\{\infty\}$,
and let $\varphi^2=dz$ represent the unique spin structure on $S^2$.
Let $P=[a_1] + \dots + [a_{n-1}] + [\infty]$ with the $a_i\in\bbC$ distinct.
We have $\calH=0$
since there are no holomorphic differentials on the sphere.
A basis for $\calF$ is 
\[
 \{\upto{t}{1}{n-1},t_n\} = \left\{\frac{\varphi}{z-a_1},\dots,
    \frac{\varphi}{z-a_{n-1}},\varphi\right\}.
\]
These sections are in $\calF$ since 
\[
(t_n) = -[\infty],\spaceout (t_i) = -[a_i],
\]
and are independent because they have distinct poles, and so are a basis for
$\calF$ since $\dim\calF=n$.
By the local calculation (\ref{local}) for $\skewform$,
\[
\skewform(t_i,t_j) = \left\{
\begin{array}{ll}
\arraystrut\displaystyle
\frac{1}{a_j-a_i} & (1\le i \le n-1; \;\;  1\le j \le n-1;\;\;i\ne j),\\
-1 & (1\le i\le n-1;\;\;j=n),\\
1 & (i=n;\;\;1\le j\le n-1),\\
0 & (i=j).
\end{array}
\right.
\]

\section{Genus zero surfaces with four or six embedded planar 
ends}\label{4ends}

The family of minimal genus zero surfaces with four embedded
planar ends was computed first in \cite{Bryant1}.
A different computation is included here for completeness, as
well as an explicit computation of the family of such surfaces
with six ends.

\begin{theorem}
The space $\Sigma_4$ of complete minimal immersions of spheres punctured
at four points into $\bbR^3$ with finite total curvature and embedded
planar ends is $S^1\times H^3$.
\end{theorem}

\begin{pf}
Let $z$ be the standard conformal coordinate on $S^2 = \bbC\cup\{\infty\}$.
By a \Moebius transformation of the Riemann sphere $S^2$,
the ends can be normalized so that
two of the ends are $0$ and $\infty$ and the product of the other two is 1.
Naming the normalized ends
\[
\{a_1=a,a_2=1/a,0,\infty\},
\]
the matrix for $\skewform$ in the basis
\[
 \left\{ \frac{1}{z-a_1}, \frac{1}{z-a_2}, \frac{1}{z}, 1 \right\}
\]
is
\[
\skewform=
\left(
\begin{array}{cccc}
\arraystrut 0 & \frac{1}{a_2-a_1} & -\frac{1}{a_1} & -1\\
\arraystrut \frac{1}{a_1-a_2} & 0 & -\frac{1}{a_2} & -1\\
\arraystrut \frac{1}{a_1} & \frac{1}{a_2} & 0 & -1\\
\arraystrut 1 & 1 & 1 & 0
\end{array}
\right)
\]
(see Section \ref{compute}).
The pfaffian (see Appendix \ref{pfaffian}) of $\skewform$ computes to
a nonzero multiple of 
\[
(a^2 - \sqrt{3} a + 1)(a^2 + \sqrt{3} a + 1).
\]
This pfaffian must be zero in order for $\ker\skewform=\calK$ to be
at least two-dimensional and hence to generate surfaces.
Setting this pfaffian to zero yields four interchangeable solutions for $a$,
one of which is
\[
a = (\sqrt{3}+i)/2.
\]
With $\varphi^2=dz$ as usual,
a basis for $\calK$ is
\[
t_1 = \left(\frac{\sqrt{3} z - 1}{z(z^2-\sqrt{3} z + 1)}\right)\varphi 
\andspace
  t_2 = \left(\frac{z(z-\sqrt{3})}{z^2-\sqrt{3} z + 1}\right)\varphi.
\]
Thus there is a single $S^1\times H^3$ family of immersions, as in Theorem \ref{theorem: moduli space}, and these have no periods since we are on $S^2$.
\end{pf}
When there are six ends, the conformal type of the domain is no longer unique:
\label{6ends}
\begin{theorem}
The space $\Sigma_6$ of complete  minimal immersions of spheres punctured 
at six points into space with finite total curvature and embedded planar ends 
is $V\times \left( S^1\times H^3\right)$, where $V$ is a complex algebraic surface.
\end{theorem}

\begin{pf}
On the sphere $S^2 = \bbC\cup\{\infty\}$ with standard conformal
coordinate $z$, the ends can be normalized so that two of the ends
are at 0 and $\infty$, and the product of the remaining four ends is 1.
With this normalization, let the ends be $\{a_1,a_2,a_3,a_4,0,\infty\}$.
Set
\[
\begin{array}{rcl}
\sig_1 \!\!\andeq \,\;\;a_1+a_2+a_3+a_4,\\
\sig_2 \!\!\andeq \!\!\!-(a_1 a_2 + a_1 a_3 + a_1 a_4 + a_2 a_3 + a_2 a_4 + 
a_3 a_4),\\
\sig_3 \!\!\andeq \,\;\;a_1 a_2 a_3 + a_1 a_2 a_4 + a_1 a_3 a_4 + a_2 a_3 a_4.
\end{array}
\]

The pfaffian (see Appendix \ref{pfaffian}) of $\skewform$ is
\begin{equation}\label{det6}
\pfaffian\skewform = \tau_1\tau_3 + \sig_1\sig_3 - 20,
\end{equation}
where
\[
\tau_1 = \sig_1^2 + 3 \sig_2\andspace \tau_3 = \sig_3^2 + 3 \sig_2.
\]
The condition that the pfaffian be 0 
defines the algebraic subvariety 
\[
V = \{(\sig_1,\sig_2,\sig_3)\subset(\bbC{\bbP}^1)^3\suchthat
    \pfaffian\skewform = 0\}
\]
of codimension 1.  Each point $(\sig_1,\sig_2,\sig_3)$ of $V$ yields a basis
\begin{eqnarray*}
t_1 \andeq \left(\frac{b_3 z^3+b_2z^2+b_1z+b_0}
     {z(z^4-\sig_1z^3-\sig_2z^2-\sig_3z+1)}\right)\varphi,\\
t_2 \andeq \left(\frac{z(c_3 z^3+c_2z^2+c_1z+c_0}
     {z^4-\sig_1z^3-\sig_2z^2-\sig_3z+1}\right)\varphi,
\end{eqnarray*}
for $\calK_{(\sig_1,\sig_2,\sig_3)}$ where 
\[
\begin{array}{l}
b_0 = \sig_2,\\ 
b_1 = -\sig_2 \sig_3,\\ 
b_2 = \sig_2\tau_3 - 2 \sig_1\sig_3 - 10,\\
b_3 = \sig_1\tau_3 + 5 \sig_3,
\end{array}
\spaceout
\begin{array}{l}
c_0 =  \sig_3\tau_1 + 5 \sig_1,\\
c_1 = \sig_2\tau_1 - 2 \sig_1\sig_3-10,\\
c_2 = -\sig_1 \sig_2,\\
c_3 = \sig_2,
\end{array}
\]
and where $\varphi^2=dz$.  The $S^1\times H^3$ family arises as before.
\end{pf}

That the four- and six-ended families are immersed follows from
Lemma~\ref{lemma:sphere} below, which in turn follows directly
from the definitions of the spaces in equation (\ref{setup}).

\begin{lemma}\label{lemma:sphere}
On the sphere with its unique spin structure $S$,
let $P_1 = \sum [p_i]$ as in equation (\ref{setup}), and
{\em $P_2 = P_1 + [a],$ $(a\not\in\supp (P_1))$}.  Let
$\calF_i = \calF_{S^1,P_i,S}$ and $\calK_i = \calK_{S^2,S,P_i}$ $(i=1,2)$
as in equation (14.11).  Then
$\calK_2\cap\calF_1 = \{s\in\calK_1\suchthat s(a)=0\}$.
\end{lemma}

Now, to see why this implies the above examples are immersed, let
$P_1$ be the divisor of ends of even degree $n<9$, and let
$(s_1, s_2)$ be the spinor representation of a minimal branched
immersion.  Supposing this surface is not immersed, let
$a$ be a branch point of the surface, and
set $P_2 = P_1+[a]$.
Then $s_1$ and $s_2$ are independent sections in $\calK_1$ 
and $s_1(a)=0$, $s_2(a)=0$, so
by Lemma~\ref{lemma:sphere} above, $s_1$, $s_2\in\calK_2$.  
Applying Lemma~\ref{dimension lemma}(ii),  we have that 
\[
2\le\dim\calK_2\le[\sqrt{n}]\le2,
\]
so $\dim\calK_2=2$.  This contradicts the fact that $n+1-\dim\calK_2$ is 
even (Lemma~\ref{dimension lemma}(i)).

\section{Projective planes with three embedded planar 
ends}\label{projective}

It was shown in \cite{Kusner2} that any finite-total-curvature  minimal immersion of a
punctured real projective  plane with embedded ends has only planar ends,
and has at least three of them.  Hence those which are the subject of
the following theorem are the examples of minimal projective planes
with the fewest number of embedded ends.  One method for determining
the moduli space of 
minimally immersed
projective planes punctured at three points was given in
\cite{Bryant2}.  Here we provide another description of this moduli
space using the spinor representation.  Note that all these surfaces
compactify to give surfaces minimizing $W = \int H^2 dA$ among all
immersed real projective planes
\cite{Kusner1}, with minimum energy $W=12\pi$.

\begin{theorem}\label{theorem:projective}
Let $\Pi_3\subset\Sigma_6$ be the moduli space of complete 
minimal immersions of real projective planes 
punctured at three points
with finite total curvature and embedded planar ends modulo Euclidean
similarities.  Then
\begin{enumerate}
\item $\Pi_3$ is homeomorphic to a closed disk with one
point $M_0$ removed from the boundary;
\item the point $M_0$ represents the \Moebius strip with total curvature $-6\pi$
in the sense that if $\gam:\bbR^+\to\Pi_3$ is a curve with 
$\lim_{t\rightarrow\infty} \gam(t) = M_0$, then there is a one-parameter
family of immersions $X_t$ parametrizing the surfaces $\gam(t)$ such that as
$t\rightarrow\infty$, $X_t$ converges uniformly on compact sets to
a parametrization of the \Moebius
strip;
\item the surfaces with non-trivial symmetry groups are represented by
the boundary of the disk, which represents
a one-parameter family of surfaces which have a
line of reflective symmetry;  among these, the only surfaces 
with larger symmetry groups
(other than $M_0$) are two surfaces which have, respectively,
the symmetry groups $\bbZ_2\times\bbZ_2$, and $D_3$, the dihedral 
group of order $6$.
\end{enumerate}
\end{theorem}

\begin{pf*}
{\it Proof of (i)}
The two-sheeted covering of the projective plane is 
the Riemann sphere $S^2=\bbC\cup\{\infty\}$, with
order-two orientation-reversing deck transformation $I(z)=-1/\ov{z}$.
By a motion in $\PSU(2)$ the six preimages on the sphere of three points in
the projective plane can be normalized as in Section \ref{6ends} to be
\[
\{a_1,I(a_1),a_2,I(a_2),0,\infty\}
\] 
with the product of the first four equal to $1$. With this choice, following
the notation of Section \ref{6ends}, we have
\[
\sig_2\in\bbR;\spaceout \sig_3 = -\ov{\sig}_1; \spaceout 
  \tau_3 = \ov{\tau}_1.
\]

For each choice of ends satisfying equation (\ref{det6}), up to 
dilations and isometries of space there is a unique minimal immersion
of the projective plane, whose spinor representation is given by
$\sqrt{i}(t_1,t_2)$, with $t_1$, $t_2$ as in Section \ref{6ends}.
For if $\sqrt{i}(\hat{t}_1,\hat{t}_2)$ is the spinor representation of
another immersion with the same ends, then a motion in
$\bbC^*\times\PSL(2,\bbC)$ can make $\hat{t}_1=t_1$, and the compatibility
condition in Theorem~\ref{nonorientable theorem}
forces $\hat{t}_2 = \pm t_1$.
Hence the moduli space $\Pi_3$ can be parametrized as a quotient
space of
\[
\Gamma = \{(\sig_1,\sig_2)\in\bbC\times\bbR\suchthat 
\tau_1\tau_3 + \sig_1\sig_3 - 20=0,\;\;\sig_3 = -\ov{\sig}_1\},
\]
where $\sig_1$, $\sig_2$, $\sig_3$ are the symmetric polynomials of
the ends defined in Section \ref{6ends}.  The desired moduli space is a 
quotient space of $\Gamma$, since permutations of the ends give rise to
the same surface.

Since the parameters $\sig_1$ and $\sig_2$ depend on the particular 
normalization of the ends made in Section \ref{6ends}, new parameters
are chosen, namely the three direction cosines $(c_1,\,c_2,\,c_3)$
of the angles between the ends $0$, $a_1$ and $a_2$,
viewed as vectors in $S^2 \subset \bbR^3$.
With these parameters, the determinant of $\skewform$ becomes,
up to a non-zero multiple,
\[
(c_1^2 + 3)(c_2^2 + 3)(c_3^2 + 3)-32(c_1 c_2 c_3 + 1).
\]
The surface
\[
 \Gamma=\left\{(c_1,c_2,c_3)\in \bbR^3 \left|\right.
 (c_1^2 + 3)(c_2^2 + 3)(c_3^2 +
3)-32(c_1 c_2 c_3 + 1) = 0\right\}
\]
in the cube 
\[
C = \left\{(x,y,z)\in \bbR^3 | -1 < x,y,z < 1 \right\}
\]
is a tetrahedron-like object
but with smoothed edges and (omitted) vertices at $(1,1,1)$, $(1,-1,-1)$,
$(-1,1,-1)$, and $(-1,-1,1)$.  

The moduli space $\Pi_3$ is diffeomorphic to a quotient of 
$\Gamma$ which arises from permutations of the ends.
A choice $c = (c_1,c_2,c_3)$ determines a set of six ends on the
double-covering sphere.
The group of rotations of the cube is the
order-24 permutation group $S_4$ generated by two kinds of elements:
\begin{enumerate}
\item[$\bullet$] permuting the three numbers $(c_1,c_2,c_3)$,
\item[$\bullet$] negating any two of the three numbers $(c_1,c_2,c_3)$.
\end{enumerate}
Action under this group determines the same six ends. 
Hence $D = \Gamma/S_4$ is a representation of the moduli space $\Pi_3$ of minimal
projective planes with three embedded planar ends.
%
$D$ can be shown to be topologically a
closed disk with the point corresponding to the corner $(1,1,1)$ of the
cube removed.

{\it Proof of (ii).}
The minimal \Moebius strip with total curvature $-6\pi$,
found in \cite{Meeks1}, has spinor representation
\[
G(w)\sqrt{dw} = \sqrt{i}(-(w+1)/w^2,w-1)\sqrt{dw}
\]

Let $(\sig_1(s),\sig_2(s)):\bbR^+\to\Gamma$ be a proper curve.  
It follows from the reality of $\sig_2$ that
\[
\lim_{s\rightarrow \infty} \frac{1}{\sig_1(s)} =
  \lim_{s\rightarrow \infty} \frac{1}{\sig_2(s)} =
  \lim_{s\rightarrow \infty} \frac{\sig_1(s)}{\sig_2(s)} = 0,
\]
and by a permutation of the ends we can assume
\[
\lim_{s\rightarrow \infty} \frac{\ov{\sig_1(s)}}{\sig_1(s)} = 1.
\]
Further,
\[
\lim_{s\rightarrow\infty}\left|\frac{\tau_1(s)}{\sig_1(s)}\right| = 1,
\]
since 
\[
\left|\frac{\tau_1}{\sig_1}\right|^2 = 
   -\frac{\tau_1 \tau_3}{\sig_1 \sig_3} = 
  1-\frac{20}{|\sig_1^2|}.
\]
Now choose a function $\al:\bbR^+\to S^1\subset\bbC$ such that
\[
\lim_{s\rightarrow\infty}\left(\frac{\tau_1(s)}{\sig_1(s)} - 
   \al(s)\right) = 0,
\]
and so
\[
\lim_{s\rightarrow\infty}\left(\frac{\tau_3(s)}{\sig_1(s)} - 
   \ov{\al(s)}\right) = 0.
\]

Let $X$ be defined by 
\[
X(z)\sqrt{dz} = \frac{\sqrt{i}}{\sig_1}(t_1,t_2),
\]
where $t_1$, $t_2$ are as in Section \ref{6ends}.
A careful reparametrization and rotation of the surface generated by
$X(z)\sqrt{dz}$ converges uniformly in compact sets to the \Moebius strip
given above:  Let $z=\al w$, and
\[
A_\al = 
\left(\begin{array}{cc} a^{3/2} & 0\\ 0 & \al^{-3/2}\end{array}\right).
\]
Then
\[
A_{\al} X(z) \sqrt{dz} = A_{\al} \sqrt{\al} X(\al w)\sqrt{dw}
\]
is the appropriate reparametrization and rotation.  This amounts to showing
\[
\lim_{s\rightarrow\infty} A_{\al(s)}\sqrt{\al(s)} X(\al(s) w) = G(w)
\]
uniformly in compact sets not containing the ends, which
follows by a calculation using the limits above.

{\it Proof of (iii).}
To find the surfaces in $\Pi_3$ which have non-trivial symmetry groups
as surfaces in space,
let $G=\bbZ_2\times\PSU(2)\cong\Orth(3)$ be the group
of conformal and anti-conformal diffeomorphisms of $\bbC\cup\{\infty\}=S^2$
with the property that any $\xi\in G$ commutes with $I$.  Via stereographic
projection, $G$ can be thought of as the isometry group of $S^2\subset
\bbR^3$,
so $\xi\in G$ satisfies $a\cdot b = \xi a \cdot \xi b.$
The group of symmetries of the minimal surface in space induces a 
subgroup $H\subset G$ acting on the domain $S^2$. 
Moreover,
the subgroup $H\subset G$ which permutes the ends is isomorphic to the subgroup
$K\subseteq S_4$ which fixes the point $(c_1,c_2,c_3)$ representing the
ends, since $\xi\in H$ preserves the inner product defining the cosines
$c_1$, $c_2$, $c_3$.  

The point of all this is that
the symmetry group of a surface represented by 
$(c_1,c_2,c_3)\in\Pi_3$ can be determined by finding the subgroup of $S_4$ 
which fixes $(c_1,c_2,c_3)$.  Using this method, the symmetric surfaces other than
the \Moebius strip at $(1,1,1)$ are
\begin{enumerate}
\item[$\bullet$] elements of $\partial D$, each with a line of reflective symmetry,
\item[$\bullet$] 
$(\sqrt{5}/3,0,0)\in\partial D$ with symmetry group $\bbZ_2\times\bbZ_2$,
\item[$\bullet$] $(c,c,-c)\in\partial D$ with symmetry group $S_3 = D_3$ 
\end{enumerate}
The last (and most symmetric) of these is a surface described in
\cite{Kusner2}.
\end{pf*}

\section{$\skewform$ on the twisted torus}
\label{du}

For the non-example in Section \ref{3ends}, and for the example in
Section \ref{4ends2}, it is necessary to compute a basis
for $\calF$ for the twisted torus (see Section \ref{spin on tori}), 
and the matrix for $\skewform$ in this basis.
On the torus $\bbC/\{2\om_1,2\om_3\}$ with the standard coordinate $u$, let
$S$ be the spin structure corresponding to the twisted torus,
that is, represented by the holomorphic differential
$\varphi_0^2=du$. Let $P =[a_1]+\dots+[a_n]$ and set $\om_2=\om_1 + \om_3$
throughout the remainder of Part III.

To show that
$\calH = \{c\varphi_0\suchthat c\in\bbC\}$, let $t\in\calH$.
Then $0\le (t) = (t/\varphi_0) + (\varphi_0) = (t/\varphi_0).$
Hence $t/\varphi_0$ is a holomorphic function on the torus, so it is constant.
A basis for $\calF$ is $\{t_0,\upto{t}{1}{n-1}\}$, where
\begin{eqnarray*}
t_0 \andeq \varphi_0,\\
t_i \andeq\left( \zeta(u-a_i) - \zeta(u) + \zeta(a_i)\right)\varphi_0,\\
\andeq \frac{1}{2}
   \left(\frac{\wp'(u)+\wp'(a_i)}{\wp(u)-\wp(a_i)}\right)\varphi_0
\end{eqnarray*}
(see equation (\ref{zeta})).  These are in $\calF$ because
\begin{eqnarray*}
(t_0) \andeq 0\ge -P,\\
(t_i) \andeq [x_i] + [y_i] - [a_i] - [0]\ge -P
\end{eqnarray*}
where $x_i$ and $y_i$ are the zeros of $\wp'(u)+\wp'(a_i)$ other than $-a_i$.
These Sections are independent because they have distinct poles, and
hence span $\calF$  since $\dim \calF = n$.
To compute $\skewform$ in this basis,
first compute the expansions of $t_i$ at $\upto{a}{0}{n-1}$ 
(assume $i$, $j\ne0$):
\[
\begin{array}{rcl}
t_i \!\!\!\andeq\!\!\! (-u^{-1} + \ooo(u))\varphi_0,\\
t_i \!\!\!\andeq\!\!\! ((t_i/\varphi_0)(a_j) +\ooo(u))\varphi_0
    \spaceout (i\ne j),\\
t_i \!\!\!\andeq\!\!\! (u-a_i)^{-1}\varphi_0.
\end{array}
\]
Using equation (\ref{local}), we have 
\[
\skewform(t_i,t_j) = \left\{
\begin{array}{ll}
\arraystrut\displaystyle
\left.\frac{t_i}{\varphi_0}\right|_{a_j} &
     \hbox{($i\ne0$; $j\ne 0$; $i\ne j$),}\\
0 & \hbox{(otherwise).}
\end{array}
\right.
\]
\section{$\skewform$ on the untwisted tori}
\label{nondu}

As above, it is also necessary to exhibit a basis for
$\calF$ on the untwisted tori (see Section \ref{spin on tori}),
as well as the matrix for $\skewform$ in this basis.
On the torus $\bbC/\{2\om_1,2\om_3\}$ with the standard conformal coordinate
$u$,
fix $r\in\{1,2,3\}$ and choose the spin structure on the untwisted torus,
represented by
\[
\varphi_r^2 = \frac{du}{\wp_r(u)},
\] where
$\wp_r(u) = \wp(u) -  \wp(\om_r)$.
Let $P=\sum[a_i]$ with the $a_i\in T\setminus\{0,\om_r\}$ distinct.

For this choice of spin structure, $\calH=0$.  To show this,
first note first that $(\varphi_r) = [0]-[\om_r]$.
If $t\in\calH$, then
\[
0 \le (t) = (t/\varphi_r) + (\varphi_r) = (t/\varphi_r) + [0] - [\om_r].
\]
It follows that $(t/\varphi_r) \ge [\om_r]-[0]$.
But since $t/\varphi_r$ is a function, the degree of its divisor is 0.
Hence $(t/\varphi_r) = [\om_r]-[0]$.  But this is impossible by Abel's
theorem on the torus: for an elliptic function $f$,
if $(f) = \sum n_i[p_i]$ (as a formal sum) then
$\sum n_i p_i = 0$ (as a sum in $\bbC$).

A basis for $\calF$ is $\{\upto{t}{1}{n}\}$, where
\begin{eqnarray*}
t_i(u) \andeq\left( \zeta(u-a_i) - \zeta(u) -\zeta(\om_r-a_i)+\zeta(\om_r)
  \right)\varphi_r\\
\andeq \frac{1}{2}\left(
\frac{\wp_r(u)\wp_r'(a_i) + \wp_r'(u)\wp_r(a_i)}
    {\wp_r(a_i)(\wp_r(u)-\wp_r(a_i))}\right)\varphi_r
\end{eqnarray*}
(see equation (\ref{zeta})).  These are in $\calF$ because 
$(\varphi_r) = [0]-[\om_r]$, so
$(t_i) = [a_i-\om_r]-[a_i]\ge -P,$
and are independent because their poles are distinct, so
they span $\calF$  since $\dim \calF = n$.
The expansions of $t_i$ at $\upto{a}{1}{n}$ are
\[
\begin{array}{rcl}
t_i \!\!\andeq\!\!\! ((t_i/\varphi_r)(a_j) + \ooo(u-a_j))\varphi_r
   \spaceout (i\ne j),\\
t_i \!\!\!\andeq\!\!\! ((u-a_i)^{-1} + \ooo(u-a_i))\varphi_r.
\end{array}
\]
Using the local expression (\ref{local}) for $\skewform$, we have 
\[
\skewform(t_i,t_j) = \left\{
\begin{array}{ll}
\arraystrut\displaystyle
\left.\frac{t_i}{\varphi_r}\right|_{a_j} & (i\ne j),\\
0 & (i=j).
\end{array}
\right.
\]

A particularly simple situation arises when the ends come in 
pairs $a$ and $-a$.  Assume $n=2m$ and  $a_{m+i} =
-a_i$ $(i = 1,\dots,m)$.  In this case, a simpler basis is 
$\{\upto{\hat{t}}{1}{m},\upto{\hat{t}}{m+1}{2m}\}$, where for $1\le i\le m$,
\[
\begin{array}{rcrcl}
\arraystrut
\hat{t}_i \!\!\!\andeq\!\!\! \disp
    \frac{\wp_r(a_i)}{\wp_r'(a_i)}\left(t_i-t_{m+i}\right)\varphi_r 
	\!\!\!\andeq\!\!\! 
    \disp\left(\frac{\wp_r(u)}{\wp_r(u)-\wp_r(a_i)}\right)\varphi_r,\\
\arraystrut
\hat{t}_{m+i}\! \!\!\andeq\!\!\! 
    (t_i+t_{m+i})\varphi_r \!\!\!\andeq\!\!\! 
    \disp\left(\frac{\wp_r'(u)}{\wp_r(u)-\wp_r(a_i)}\right)\varphi_r.
\end{array}
\]
In this basis, the matrix for $\skewform$ becomes
\[
\left(
\begin{array}{c|c}
0 & W\\ \hline
\arraystruti
-\trans{W} & 0
\end{array}
\right),
\]
where $W$ is given by
\[
W_{ij} = \left\{
\begin{array}{ll}
\arraystrut
\disp \frac{4}{\wp_r(a_i)-\wp_r(a_j)} & (i<j),\\
\arraystrut
\disp \frac{4}{\wp_r(a_j)-\wp_r(a_i)} & (i>j),\\
\arraystrut
\disp \frac{\wp_r(a_i)^2-c_p c_q}{\wp_r(a_i)(\wp_r(a_i)-c_p)(\wp_r(a_i)-c_q)}
    & (i=j)
\end{array}
\right.
\]
and $c_p= e_p-e_r$, $c_q = e_q - e_r$, $\{p,q,r\} = \{1,2,3\}$.
Note that the entries of $W$ are entirely free of 
$\wp_r'$.

A useful property of the basis above is as follows:
let $L:M\longrightarrow M$ be defined as $L(u)=-u$;
then for $i\le m$ and $j\ge m+1$, 
\[
 L^*(\hat{t}_i \hat{t}_j ) = \hat{t}_i \hat{t}_j ,
\]
so
\[
\int_{\gam_k} \hat{t}_i \hat{t}_j  = \int_{\gam_k} L^*(\hat{t}_i \hat{t}_j ) = 
\int_{L(\gam_k)} \hat{t}_i \hat{t}_j  = -\int_{\gam_k} \hat{t}_i \hat{t}_j,
\]
and so
\[
\int_{\gam_k} \hat{t}_i \hat{t}_j = 0\spaceout(i\le m;\ j\ge m+1;\ k=1,3).
\]

\section{Non-existence of tori with three planar ends}
\label{3ends}

An outline of the proof of the non-existence of three-ended tori,
twisted or untwisted, is given.

\begin{theorem}\label{theorem:3ends}
There does not exist a complete minimal branched immersion 
of a torus into space with finite total curvature
and three embedded planar ends.
\end{theorem}

\begin{pf*}{\it Sketch of proof}
The proof is divided into two cases: for the twisted torus
there exist immersions with periods, but the periods cannot
be made purely imaginary;
for the untwisted torus, there are not even periodic examples.

First consider the more difficult case of the twisted torus.
With everything as in Section \ref{du},
let $\{0,a_1,a_2\}$ be the set of
ends, and let $p_i=\wp(a_i)$, $p_i'=\wp'(a_i)$.
The condition $\dim\calK\ge 2$ puts the following
condition on the placement of the ends:
\[
g_2 = 4(p_1^2 + p_1 p_2 + p_2^2),
\]
where $g_2$ is the constant in the differential equation
$(\wp')^2 = 4 \wp^3 - g_2 \wp - g_3$.
To see this, first note that $\ker\skewform = \calK\oplus\calH$ and
$\dim\calH = 1$. Hence  if $\dim\calK=2$ then $\skewform\equiv 0$.
Assume first that $a_1+a_2\ne 0$. Then $p_1-p_2\ne0$, and
the entries of $\skewform$
indicate that $p_1'+p_2'=0$
Hence
\[
(p_1')^2 = 4 p_1^3 - g_2 p_1 - g_3
\]
and
\[
(p_2')^2 = 4 p_2^3 - g_2 p_2 - g_3
\]
are equal, and the desired condition follows.
The condition also obtains in the case that $a_1+a_2=0$;
this can be shown as a limiting case of the above.

Changing basis now to simplify the period equations,
let 
\begin{eqnarray*}
\hat{t}_1 \andeq t_1 + \ep t_2,\\
\hat{t}_2 \andeq t_1 + \ep^2 t_2,
\end{eqnarray*}
where $\ep = (-1+\sqrt{3})/2$.
With $\gam_1$, $\gam_3$ the closed curves parallel to
$\om_1$, $\om_3$ respectively (as in Theorem~\ref{kleintype}),
the integrals relevant to the periods are (for $k=1$, $3$)
\[
\int_{\gam_k} \hat{t}_1^2 = -6 q_1 \om_k,\spaceout
\int_{\gam_k} \hat{t}_1\hat{t}_2 = -6 \eta_k,\spaceout
\int_{\gam_k} \hat{t}_2^2 = -6 q_2 \om_k,
\]
where
\begin{eqnarray*}
q_1 \andeq -( (\ep-\ep^2)p_1 + (\ep-1)p_2)/3,\\
q_2 \andeq -( (\ep^2-\ep)p_1 + (\ep^2-1)p_2)/3,\\
q_1 q_2 \andeq (p_1^2 + p_1 p_2 + p_2^2)/3 = g_2/12.
\end{eqnarray*}
A choice of a pair of independent Sections from $\calK$ can be
normalized by the action of $\bbR^*\times\SU(2)$ to be
\begin{eqnarray*}
s_1 \andeq z_1 \hat{t}_1 + \hat{t}_2,\\
s_2 \andeq z_2 \hat{t}_1,
\end{eqnarray*}
with $z_1$, $z_2\in\bbC$.
Then the period equations (\ref{period1}) can be written
\begin{eqnarray*}
\left(\begin{array}{c} 2 z_1\\ z_1^2 q_1 + q_2\end{array}\right) - B
\left(\begin{array}{c} 0\\ \ov{q}_1\ov{z}_2^2\end{array}\right) \andeq 0,\\
\left(\begin{array}{c} z_2\\q_1 z_1 z_2\end{array}\right) + B
\left(\begin{array}{c} \ov{z}_2\\ \ov{q}_1\ov{z}_1\ov{z}_2^2
     \end{array}\right) \andeq 0,
\end{eqnarray*}
where
\[
B = A^{-1}\ov{A} =
\left(\begin{array}{cc} a & b\\ c & d \end{array}\right);
\spaceout 
A = \left(\begin{array}{cc}\eta_1 & \om_1\\ \eta_3 & \om_3 \end{array}\right).
\]
Changing from the variables $(z_1,\, z_2)$ to $(w,\, z_2)$, this system is
equivalent to the system
\begin{eqnarray*}
w^2 + b^2 q_1 q_2 - d^2 \andeq 0,\\
2w + 2d - b^2 q_1 \ov{q}_1 \ov{z}_2^2 \andeq 0,\\
w z_2 + \ov{z}_2 \andeq 0.
\end{eqnarray*}
From these it follows that
\begin{eqnarray*}
w\ov{w}-1 \andeq 0,\\
a w^2 - \ov{a} \andeq 0,\\
-\ov{a} - a b^2 q_1 q_2 + a d^2 \andeq 0.
\end{eqnarray*}
This last condition, depending only on the conformal type of the torus and
not on $w$, $z_1$, and $z_2$,
is a degeneracy condition for the period equations.
It also follows, by an examination of the sign of $a(w -\ov{a})\in\bbR$, that
\[
 |a|>1.
\]
A delicate argument, which is omitted here,
using the expansions \cite{Lang}
\begin{eqnarray*}
g_2 \andeq \frac{\pi^4}{12 \om_1^4}\left( 1 + 240 \sum_{n=1}^{\infty}
    \sig_3(n) q^n\right),\\
\eta_1 \andeq \frac{\pi^2}{12 \om_1}\left( 1 - 24 \sum_{n=1}^{\infty}
    \sig_1(n) q^n\right),
\end{eqnarray*}
where 
\[
\sig_k(n) = \sum_{d|n}d^k; \spaceout q = e^{2 i \pi \tau};\spaceout 
\tau = \om_3/\om_1
\]
shows that the degeneracy condition is not satisfied under the constraint
$|a|>1$ over the whole moduli space of Riemann tori.
Hence no examples with three ends can be found in the case of the twisted tori.

The case of the untwisted tori is much easier.  Fix $r\in\{1,2,3\}$
and let $\varphi_r$ be as in Section \ref{nondu}.
Let $\{a_1,a_2,a_3\}$
be the ends, translated so that they avoid $\{0,\om_r\}$, and let
$\{t_1,t_2,t_3\}$ be the basis for $\calF$ given in the same section.
The condition that $\dim\calK=\dim\ker\skewform\le 2$ forces
$\skewform$ to be zero.  This means, for example, that $t_1/\varphi_r$
have zeros at $a_2$ and $a_3$.
But the zeros of $t_1/\varphi_r$ are $\om_r$ and $a_1-\om_r$, so one of
$a_2$, $a_3$ has to be $\om_r$,  contrary to the assumption.
\end{pf*}

\section{Minimal tori with four embedded planar ends}\label{4ends2}

Here the existence of families of four-ended minimal tori --- none of
which are regularly homotopic to embedded tori --- is established.
These surfaces conformally compactify to yield $W$-critical twisted
tori with $W=16\pi$ and isolated umbilics.

\begin{theorem}\label{theorem:4ends2}
For each conformal type of torus there exists a real two-dimensional
family of complete minimal immersions of the torus punctured at four points
into space with finite total curvature and embedded planar ends.  Each of
the tori is twisted.
\end{theorem}

\begin{pf}
To exhibit the family, it is first necessary to determine 
the placement of the four ends. The ends in fact must be,
up to a translation,  at the four half-lattice points.  
To show this, on the torus $\bbC/\{2\om_1,2\om_3\}$,
assume the four ends are $\{0,a_1,a_2,a_3\}$, where
$a_1$, $a_2$, $a_3$ are distinct points in the torus to be determined.
With $\varphi_0^2=du$, the matrix for $\skewform$ in the basis 
$\{t_0,t_1,t_2,t_3\} = \{\varphi_0,f_1\varphi_0,f_2\varphi_0,f_3\varphi_0\}$
of Section \ref{du} is
\[
\skewform=
\left(
\begin{array}{cccc}
0 &        0 &        0 &      0   \\
0 &        0 & f_1(a_2) & f_1(a_3)  \\
0 & f_2(a_1) &        0 & f_2(a_3)  \\
0 & f_3(a_1) & f_3(a_2) &        0
\end{array}
\right).
\]

If $\ker \skewform=\calH\oplus\calK$ is two-dimensional, then
$\dim\calK=1$, since $\dim\calH=1$, so $\calK$ is
not big enough to generate a minimal surface.
Hence to produce surfaces, $\rank \skewform$, being even, must be zero.
In this case, all the entries of the above matrix are zero;  
a look at $t_i$ shows that  $\wp'(a_i) + \wp'(a_j)=0$ for all $i\ne j$.
It follows that $\wp'(a_1) = \wp'(a_2) = \wp'(a_3) = 0$, 
so $\{a_1,a_2,a_3\} = \{\om_1,\om_2,\om_3\}$.  

With the ends fixed at
$\{0,\om_1,\om_2,\om_3\}$,
$\calF = \ker\skewform = \calH\oplus\calK$,  so 
$\{t_1,t_2,t_3\}$ is a basis for $\calK$.
The simple zeros and poles of 
$t_1$, $t_2$, and $t_3$ are shown in the following figure.
\begin{figure}[htbp]

\setlength{\unitlength}{.18in}

\centering
\begin{picture}(17,7)(0,-0.8)

\multiput(0,0)(6,0){3}{\line(1,0){4}}
\multiput(1,6)(6,0){3}{\line(1,0){4}}

\multiput(0,0)(6,0){3}{\line(1,6){1}}
\multiput(4,0)(6,0){3}{\line(1,6){1}}

\multiput(0,0)(2,0){9}{\circle*{.14}}
\multiput(.5,3)(2,0){9}{\circle*{.14}}
\multiput(1,6)(2,0){9}{\circle*{.14}}

\multiput(-.5,.4)(6,0){3}{\makebox(0,0){$\infty$}}
\put(2,.4){\makebox(0,0){$\infty$}}
\put(0,3.4){\makebox(0,0){$0$}}
\put(2.5,3.4){\makebox(0,0){$0$}}

\put(8,.4){\makebox(0,0){$0$}}
\put(6,3.4){\makebox(0,0){$0$}}
\put(8.5,3.4){\makebox(0,0){$\infty$}}

\put(14,.4){\makebox(0,0){$0$}}
\put(12,3.4){\makebox(0,0){$\infty$}}
\put(14.5,3.4){\makebox(0,0){$0$}}

\put(2,-0.7){\makebox(0,0){$t_1$}}
\put(8,-0.7){\makebox(0,0){$t_2$}}
\put(14,-0.7){\makebox(0,0){$t_3$}}

\end{picture}

\caption{Zeros and poles of $t_1$, $t_2$, and $t_3$}

\end{figure}

To solve the period problem outlined in Section \ref{periods}
it is convenient to choose a new basis
$\{\hat{t}_1,\hat{t}_2,\hat{t}_3\}$ for $\calK$
which ``diagonalizes'' the period equations.
Let
\[
\left(
\begin{array}{c} \hat{t}_1\\ \hat{t}_2\\ \hat{t}_3\end{array}\right)
=
\left(
\begin{array}{rrr}
1  & -1 & -1\\
-1 & 1  & -1\\
-1 & -1 & 1
\end{array}
\right)
\left(\begin{array}{c} t_1\\ t_2\\ t_3 \end{array} \right),
\]
or
\begin{eqnarray*}
\hat{t}_1(u) \andeq 
  (\zeta(u)+\zeta(u-\om_1)-\zeta(u-\om_2)-\zeta(u-\om_3)+
     2\zeta(\om_1))\varphi_0,\\
\hat{t}_2(u) \andeq 
  (\zeta(u)-\zeta(u-\om_1)+\zeta(u-\om_2)-\zeta(u-\om_3)+
     2\zeta(\om_2))\varphi_0,\\
\hat{t}_3(u) \andeq 
  (\zeta(u)-\zeta(u-\om_1)-\zeta(u-\om_2)+\zeta(u-\om_3)+
     2\zeta(\om_3))\varphi_0.
\end{eqnarray*}
The simple zeros and poles of $\hat{t}_1$, $\hat{t}_2$, and $\hat{t}_3$ are illustrated below.
\begin{figure}[htbp]

\setlength{\unitlength}{.18in}

\centering

\begin{picture}(17,7.5)(0,-1.2)

\multiput(0,0)(6,0){3}{\line(1,0){4}}
\multiput(1,6)(6,0){3}{\line(1,0){4}}

\multiput(0,0)(6,0){3}{\line(1,6){1}}
\multiput(4,0)(6,0){3}{\line(1,6){1}}

\multiput(0,0)(1,0){4}{\circle*{.14}}
\multiput(6,0)(1,0){4}{\circle*{.14}}
\multiput(12,0)(1,0){4}{\circle*{.14}}

\multiput(.25,1.5)(1,0){4}{\circle*{.14}}
\multiput(6.25,1.5)(1,0){4}{\circle*{.14}}
\multiput(12.25,1.5)(1,0){4}{\circle*{.14}}

\multiput(.5,3)(1,0){4}{\circle*{.14}}
\multiput(6.5,3)(1,0){4}{\circle*{.14}}
\multiput(12.5,3)(1,0){4}{\circle*{.14}}

\multiput(.75,4.5)(1,0){4}{\circle*{.14}}
\multiput(6.75,4.5)(1,0){4}{\circle*{.14}}
\multiput(12.75,4.5)(1,0){4}{\circle*{.14}}

\multiput(1,6)(1,0){4}{\circle*{.14}}
\multiput(7,6)(1,0){4}{\circle*{.14}}
\multiput(13,6)(1,0){4}{\circle*{.14}}

\multiput(-.5,-.4)(6,0){3}{\makebox(0,0){$\infty$}}
\multiput(2,-.4)(6,0){3}{\makebox(0,0){$\infty$}}
\multiput(0,3)(6,0){3}{\makebox(0,0){$\infty$}}
\multiput(2.5,3.4)(6,0){3}{\makebox(0,0){$\infty$}}

\multiput(1,.4)(2,0){2}{\makebox(0,0){$0$}}
\multiput(1.5,3.4)(2,0){2}{\makebox(0,0){$0$}}

\multiput(7.3,1.9)(2,0){2}{\makebox(0,0){$0$}}
\multiput(7.8,4.9)(2,0){2}{\makebox(0,0){$0$}}

\multiput(12.7,1.9)(1.6,0){2}{\makebox(0,0){$0$}}
\multiput(13.2,4.9)(1.6,0){2}{\makebox(0,0){$0$}}

\put(2,-1){\makebox(0,0){$\hat{t}_1$}}
\put(8,-1){\makebox(0,0){$\hat{t}_2$}}
\put(14,-1){\makebox(0,0){$\hat{t}_3$}}

\end{picture}

\caption{Zeros and poles of $\hat{t}_1$, $\hat{t}_2$, and $\hat{t}_3$}

\end{figure}

To compute the periods, use equation (\ref{standard}) to write
\begin{eqnarray*}
\hat{t}_i^2(u) \andeq 
   \left(\wp(u)+\wp(u-\om_1) + \wp(u-\om_2) + 
   \wp(u-\om_3) - 4\wp(\om_i)\right)du,\\
(\hat{t}_1\hat{t}_2)(u) \andeq 
   \left(\wp(u) - \wp(u-\om_1) - \wp(u-\om_2) + \wp(u-\om_3)\right)du,\\
(\hat{t}_1\hat{t}_3)(u) \andeq 
   \left(\wp(u) - \wp(u-\om_1) + \wp(u-\om_2) - \wp(u-\om_3)\right)du,\\
(\hat{t}_2\hat{t}_3)(u) \andeq 
   \left(\wp(u) + \wp(u-\om_1) - \wp(u-\om_2) - \wp(u-\om_3)\right)du.
\end{eqnarray*}
With $\gam_1$, $\gam_3$ the closed curves on the torus respectively
parallel to $\om_1$, $\om_3$,
the periods are
\[
P_k^{ij} = \int_{\gam_k}\hat{t}_i \hat{t}_j du = 
\left\{
\begin{array}{cc}
-8(\eta_k + \om_k e_i) & \hbox{if }i=j\\
0                      & \hbox{if }i\ne j
\end{array}
\right.
\spaceout (k=1,3),
\]
where $e_i = \wp(\om_i)$ and $\eta_k = \zeta(\om_k)$ 
(see appendix \ref{elliptic}).
In general, with
\begin{eqnarray*}
t_1 \andeq x_1 \hat{t}_1 + x_2 \hat{t}_2 + x_3 \hat{t}_3\\
t_2 \andeq y_1 \hat{t}_1 + y_2 \hat{t}_2 + y_3 \hat{t}_3
\end{eqnarray*}
the period equations (\ref{period1}) are
\begin{eqnarray*}
\sum_{1\le i,j\le 3} P_k^{ij}x_i x_j \andeq
    \ov{\sum_{1\le i,j\le 3} P_k^{ij} y_i y_j} \spaceout (k=1,3)\\
\sum_{1\le i,j\le 3} P_k^{ij}x_i y_j \andin i\bbR\spaceout (k=1,3).
\end{eqnarray*}

Now let $(i,j,k)$ be a permutation of $(1,2,3)$ and
make the particular choice
\begin{eqnarray*}
s_1 \andeq x_i \hat{t}_i + x_j \hat{t}_j,\\
s_2 \andeq \hat{t}_k.
\end{eqnarray*}
The second period equation above is satisfied for all $x_i$, $x_j$, and
the first period equation can be written in the form
\[
\left(\begin{array}{c}x_i^2\\x_j^2\end{array}\right) =
\left(\begin{array}{cc} 1 & 1\\ e_i & e_j\end{array}\right)^{-1}
B 
\left(\begin{array}{c} 1\\ \ov{e}_k \end{array}\right)
\]
where $\eta_i = \zeta(\om_i)$ and $e_i = \wp(\om_i)$ and
$B$ is defined in Section \ref{3ends}.
The condition that the surface be immersed is that $s_1$ and $s_2$ have
no common zeros.  The zeros of $s_2$ are at
$\{\om_k/2, \om_k/2+\om_1,\om_k/2+\om_2,\om_k/2+\om_3\}$,
and 
\[
\widehat{t}_m^2(\om_k/2) = \widehat{t}_m^2(\om_k/2+\om_l) = 4(e_k-e_i)
\spaceout (m=i,j;\;\;l = 1,2,3).
\]
A necessary condition that the surface branch is that
\[
 (e_k-e_i) x_i^2 - (e_k-e_j)x_j^2 = 0,
\]
or
\[
\left(\begin{array}{cc}g_2/2-3e_k^2 & -3e_k\end{array}\right)
B
\left(\begin{array}{c} 1\\ \ov{e}_k \end{array}\right)
=0.
\]
With the choice $\{i,j,k\}=\{1,2,3\}$ it can be shown that this condition is
not satisfied in the standard fundamental region of the moduli space 
of tori.  The proof uses the $q$-expansion for $g_2$ and $\eta$ given in 
Section \ref{3ends}, as well as the expansion
\[
e_1 = \frac{\pi^2}{6\om_1^2}
   \left(1 + 24\sum_{n=1}^\infty \tau(n)q^n\right),
\]
where
\[
\tau(n) = \sum_{d|n \atop d{\rm\ odd}} d.
\]

Thus we have found a single immersion of every conformal type of torus 
punctured at the half-lattice points.  Since the period conditions amount
to at most six real conditions on 12 variables, there is a real 6-parameter
family of surfaces, which modulo the action of the group in equation
(\ref{homog2})
leaves a 2-parameter family. The existence of the real two-dimensional family
follows from the fact 
that the condition of being immersed is an open analytic condition.
\end{pf}

\section{Minimal Klein bottles with embedded planar ends}\label{kleinex}

A minimal Klein bottle is constructed in this section.  Its
compactification is a $W$-critical surface with energy $W=16\pi$,
which lies in the {\em amphichiral} regular homotopy class
{\bf K}$_0\!=\,${\bf B}$\#\overline{\makebox{\bf B}}$ of Klein
bottles (cf. \cite{Kusner1}, \cite{Pinkall}).
There are no minimal Klein bottles with two embedded ends \cite{Kusner2},
and we conjecture there are none with three embedded planar ends.

\begin{theorem}\label{theorem:klein}
There exists a minimal immersion of the
Klein bottle with four embedded planar ends.
\end{theorem}

\begin{pf}
To construct this example, let 
$T = \bbC/\{2\om_1,2\om_3\}$ be a square lattice with
\[
\om_3 = i\om_1, \spaceout \om_2 = -\om_1-\om_3, \spaceout
\wp(\om_1) = 1, \spaceout \wp(\om_2) = 0, \spaceout \wp(\om_3) = -1.
\]
Let $I\colon T\to T$ be the deck transformation $I(u)=\bar{u}+\om_1$
as in Theorem~\ref{kleintype}(i) of Appendix \ref{klein}.
Let $a\in T$ be a point (yet to be determined) such that $I(a)=-a,$
and let $E=\{a_1,\ldots, a_8\}\subset T$ be the points in Table~3.
\begin{table}[htbp]
\centering
\caption{Values of $\wp$ and $\wp'$ at ends of Klein bottle}
\vspace{.1in}
$
\begin{array}{l||c|c||c}
  u             & \wp(u) & \wp'(u)  & I(u)        \\ \hline
a_1 = a         & r      & r'       & a_5 \\
a_2 = a+\om_2   & -1/r   & r'/r^2   & a_6 \\
a_3 = -ia       & -r     & -ir'     & a_4 \\
a_4 = -ia+\om_2 & 1/r    & -ir'/r^2 & a_3 \\
a_5 = -a_1      & r      & -r'      & a_1 \\
a_6 = -a_2      & -1/r   & -r'/r^2  & a_2 \\
a_7 = -a_3      & -r     & ir'      & a_8 \\
a_8 = -a_4      & 1/r    & ir'/r^2  & a_7
\end{array}
$
\end{table}

We want to construct a minimal immersion $X\colon (T\setminus E)/I
\to \bbR^3$,
\[
X(z)={\re}\int^z(s_1^2-s_2^2, i(s_1^2+s_2^2), 2s_1s_2),
\]
where $s_1,s_2$ are sections of the spin structure $S$ determined by 
$\varphi$, with
\[
\varphi^2=\frac{du}{\wp(u)-\wp(\om_2)}=\frac{du}{\wp(u)}.
\]

\paragraph{\it Step~1:\, Determination of the ends.}
Let $\{t_1,\ldots,t_8\}$,
\[
t_\alpha=\frac{\wp(u)}{\wp(u)-\wp(a_\alpha)}\varphi,\spaceout
t_{\alpha+4}=\frac{\wp'(u)}{\wp(u)-\wp(a_\alpha)}\varphi\spaceout
(1\leq \alpha\leq 4)
\]
be a basis for $\cal F$, as in Section \ref{nondu}.  The skew-symmetric
matrix for $\skewform$ in this basis is
\[
\left(
\begin{array}{c|c}
0 & W\\ \hline
\arraystruti
-\trans{W} & 0
\end{array}
\right),
\]
where $W$ is given by
\[
W=
\left(
\begin{array}{cccc}
\arraystrut
\disp\frac{r^2+1}{r(r^2-1)} & \disp\frac{4r}{r^2+1} & 
     \disp\frac{2}{r} & \disp\frac{4r}{r^2-1}\\
\arraystrut
\disp\frac{-4r}{r^2+1} & \disp\frac{r(r^2+1)}{r^2-1} & 
     \disp\frac{4r}{r^2-1} & -2r\\
\disp\frac{-2}{r} & \disp\frac{-4r}{r^2-1} & \disp\frac{-(r^2+1)}{r(r^2-1)} & 
     \disp\frac{-4r}{r^2+1}\\
\arraystrut
\disp\frac{-4r}{r^2-1} & 2r & \disp\frac{4r}{r^2+1} & 
     \disp\frac{-r(r^2+1)}{r^2-1}
\end{array}
\right).
\]

The desired sections $s_1,s_2$ lie in $\ker \skewform$, so a
necessary condition for existence is that
\[
0=\det W=\frac{(3 r^8 - 4r^6 + 50 r^4 - 4 r^2 + 3)^2}
{ (r^4-1)^2}=\frac{9(r^4+mr^2+1)^2(r^4+\ov{m}r^2+1)^2}
{(r^4-1)^2},
\]
where $m=-2(1-4\sqrt{2}i)/3$.  Let $r$ be the root of 
$r^4+mr^2+1$ in the fourth quadrant; with this choice, the
domain $T\backslash E$ is shown below.

\vspace{.12in}

\setlength{\unitlength}{0.8pt}

\begin{figure}[htbp]

\centering

\begin{picture}(160,160)(0,0)

\put(0,0){\framebox(160,160){}}
\put(0,80){\line(1,0){160}}    
\put(80,0){\line(0,1){160}}    

\put(40,10){\circle*{2.8}}      
\put(40,70){\circle*{2.8}}
\put(10,120){\circle*{2.8}}
\put(70,120){\circle*{2.8}}
\put(90,40){\circle*{2.8}}
\put(150,40){\circle*{2.8}}
\put(120,90){\circle*{2.8}}
\put(120,150){\circle*{2.8}}

\put(40,17){\makebox(0,0){$a+\om_2$}}
\put(38,63){\makebox(0,0){$-a$}}
\put(27,127){\makebox(0,0){$-ia+\om_2$}}
\put(69,110){\makebox(0,0){$ia$}}
\put(90,50){\makebox(0,0){$-ia$}}
\put(140,33){\makebox(0,0){$ia+\om_2$}}
\put(122,97){\makebox(0,0){$a$}}
\put(120,142){\makebox(0,0){$-a+\om_2$}}

\put(87,87){\makebox(0,0){$0$}}
\put(172,80){\makebox(0,0){$\om_1$}}
\put(-7,-7){\makebox(0,0){$\om_2$}}
\put(80,168){\makebox(0,0){$\om_3$}}

\end{picture}

\caption{The eight ends in the double cover of the Klein bottle}

\end{figure}
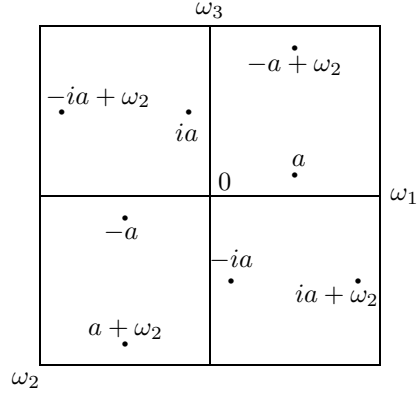

\vspace{.09in}

\paragraph{\it Step~2:\, Choosing sections $s_1,s_2$ of $S$;\,
the period equations.}
With $r$ fixed as above, rank~$\skewform$ is $4,$
and a basis for $\ker \skewform$ is $\{\hat{s}_1,\hat{s}_2,
\hat{s}_3,\hat{s}_4\}$ where
\[
\hat{s}_1 = \sum ^4_{a=1}c^\alpha_1t_\alpha, \spaceout
\hat{s}_2 = \sum ^4_{a=1}c^\alpha_2 t_\alpha, \spaceout
\hat{s}_3 = i\ov{I^*\hat{s}_1}, \spaceout
\hat{s}_4 = i\ov{I^*\hat{s}_2},
\]
\begin{eqnarray*}
c_1\andeq (2(r^2-1)^2, (r^2+1) (r^2-3), (r^2+1) (3r^2-1), -2(r^2-1)^2), \\
c_2\andeq ((r^2+1) (3r^2-1), -2(r^2-1)^2, 2(r^2-1)^2, (r^2+1) (r^2-3)),
\end{eqnarray*}
and $I^*$ is a choice of a lift of the deck transformation $I$
to the spin structure $S$.

Let
\[
\begin{array}{lr}
s_1=x_1\hat{s}_1+x_2\hat{s}_2 \\
& x_1,\,x_2\in {\bbC}\,. \\
s_2=i\ov{I^*s_1}=\ov{x}_1\hat{s}_3+\ov{x}_2
\hat{s}_4
\end{array}
\]
We want to find $x_1,\, x_2$ such that the real part of all
periods are zero.  By Theorem~\ref{kleintype}(iv) in Appendix \ref{klein}, 
the period equations reduce to the single equation
\[
0=\int_{\gam_1}s^2_1=x^2_1P^{11}_1+2x_1x_2P^{12}_1+x^2_2P^{22}_1,
\]
where
\[
P^{\alpha\beta}_k=\int_{\gam_k}\hat{s}_\alpha\hat{s}_\beta
\]
along the curve $\gam_k\colon t\longmapsto t\om_k\, (-1\leq
t\leq 1)$.
\vspace{.1in}

\paragraph{\it Step~3:\, Explicit solution of the period equation.}
The period equation above can be solved once $P^{\alpha\beta}_k$
are known.  To compute these, let
\begin{eqnarray*}
\hat{s}_1^2=\frac{1}{2}\left( -\sum_{\alpha=1}^4 A_\alpha 
\wp(u-a_\alpha) + B\right)du,
    \spaceout\textstyle  A = \sum A_\alpha\\
\hat{s}_1\hat{s}_2=\frac{1}{2}\left( -\sum_{\alpha=1}^4 C_\alpha 
\wp(u-a_\alpha) + D\right)du,
    \spaceout\textstyle  C = \sum C_\alpha
\end{eqnarray*}
as in equation (\ref{standard}). Then
\begin{eqnarray*}
P_1^{11} \andeq \int_{\gam_{1}} \hat{s}_1^2=A\eta_1 + B\om_1,\\
P_1^{12} \andeq \int_{\gam_{1}} \hat{s}_1\hat{s}_2=C\eta_1 + D\om_1,\\
P_3^{11} \andeq A\eta_3 + B\om_3 = i(-A\eta_1 + B\om_1) \\
P_3^{12} \andeq C\eta_3+D\om_3=i(-C\eta_1+D\om_1) \\
\eta_k \andeq -\frac{1}{2}\int_{\gam_{k}}\wp(u)du.
\end{eqnarray*}

Let $J\colon T\rightarrow T$ be defined by $J(u)=iu$, and let
$J^*$ be a lift of $J$ to $S$.  Then
\[
\hat{s}_1=\sqrt{i}J^*\hat{s}_2, \quad \hat{s}_2=\sqrt{i}J^*
\hat{s}_1
\]
for some choice of $\sqrt{i}$.  Then
\[
P^{12}_1=\int_{\gam_{1}}\hat{s}_1\hat{s}_2=\int_{\gam_{1}}iJ^*
\hat{s}_1\hat{s}_2=i\int_{J_{(\gam_1)}}J^*\hat{s}_1\hat{s}_2=i
\int_{\gam_{3}}\hat{s}_1\hat{s}_2=iP_3^{12},
\]
so $D=0$. Again,
\[
P_1^{22}=\int_{\gam_{1}}\hat{s}^2_2=\int_{\gam_{1}}iJ^*\hat{s}_1
^2=i\int_{J_{(\gam_{1})}}\hat{s}_1^2=i\int_{\gam_{3}}\hat{s}_1^2=iP_3^{11},
\]
so $P_1^{22}=A\eta_1-B\om_1$.

Having computed $P_1^{11},P_1^{12},P_1^{22}$ in terms of $ A,B,C $, we
compute $ A,B,C $ by expanding $\hat{s}_\alpha\hat{s}_\beta/du$ in two ways
and equating coefficients.  On the one hand, by the definition of
$\hat{s}_\alpha$, we have
\[
\hat{s}_a\hat{s}_\beta/du=\sum_{\gam,\delta}\frac{c^\gam_\alpha
c^\delta_\beta\wp(u)}{(\wp(u)-\wp(a_\gam))(\wp(u)-\wp(a_\delta))}
\; (1\leq \alpha,\,\beta\leq 2;\,\, 1\leq \gam,\,\delta\leq 4).
\]
Using the formula (for $\wp'(u_0)$ finite and non-zero)
\[
\frac{1}{\wp(u)-\wp(u_0)}=\frac{1/\wp'(u_0)}{u-u_0}+\cdots,
\]
we get the expansion at $a_\gam$
\[
\hat{s}_a\hat{s}_\beta/du=\frac{c^\gam_ac^\gam_\beta\wp(a_\gam)/
(\wp'(a_\gam))^2}{(u-a_\gam)^2}.
\]
On the other hand we have the expansions at $a_\gam$
\[
\hat{s}^2_1/du = \frac{-A_\gam/2}{(u-a_\gam)^2}
\andspace
\hat{s}_1\hat{s}_2/du = \frac{-C_\gam/2}{(u-a_\gam)^2}.
\]
Equating coefficients,
\begin{eqnarray*}
A_\gam \andeq -2(c_1^\gam)^2 \wp(a_\gam)/(\wp'(a_\gam))^2 \\
C_\gam \andeq -2 c_1^\gam c_2^\gam \wp(a_\gam)/(\wp'(a_\gam))^2,
\end{eqnarray*}
so
\begin{eqnarray*}
A \andeq \sum A_i = -32 r^2(r^4+4r^2+1)/3 \\
C \andeq \sum C_i = -2(r^4-1)^2.
\end{eqnarray*}
To compute $B$, note that $s_1$ has a zero at $0$ to get
\[
B=\sum A_\gam \wp(a_\gam)=4r(r^2+1)^3.
\]
This solves the period equation.

Finally, that the immersion is unbranched is the condition that
$s_1,\,s_2$ have no common zeros.
This amounts to the condition that if $u_0$ is a zero of $s_1$, then
$I(u_0)$ is not.
By using the identity
\[
\ov{I^*\wp} = \frac{\wp+1}{\wp-1},
\]
this can be checked by
setting $s_1$ to zero, and solving numerically
the cubic in $\wp$ which results.
\end{pf}
\nopagebreak
%
%

\typeout{_______________________________________________ appendix.tex}

\appendix   

\section{Winding numbers and quadratic forms}\label{quadratic appendix}

In this appendix, we sketch the proof that $q_S(c) = w(\al,\,v)+1$
defines a $\bbZ_2$-quadratic form associated to the spin structure $S$ (Theorem
\ref{quadratic theorem}).

\noindent
{\em Proof.}  Let $\al_0$, and  $\al_1:S^1\to M$ be embedded representatives
(see \cite{Meeks2}, \cite{Kauffman2}) of $c\in
H_1(M,\bbZ_2)$.  Let $v_0$, $v_1$ be smooth nonzero vector fields which
lift along $\al_0$, $\al_1$ respectively to sections of the spin
structure $S$.  Let $\al_t$ $(t\in[0,1])$ be a homotopy of $\al_0$ and
$\al_1$.  Extend $v_0$ to a smooth nonzero vector field $v$ in an annulus
containing the image of $\al_t$.
Then $w(\al_t,v)$ is a continuous function of $t$, and an integer, hence
it is constant.  In particular,
\[
w(\al_0,v_0) = w(\al_1,v).
\]
But $v=v_0$ lifts along $\al_0$ to a smooth section of $S$.
So $v$ must also lift along $\al_1$.
Since $v_1$ also lifts along $\al_1$, by the lemma below
\[
w(\al_1,v) = w(\al_1,v_1).
\]
Thus
\[
w(\al_0,v_0) = w(\al_1,v_1),
\]
showing that $q_S$ is well-defined.

Now, to show $q_S$ is quadratic, let $\al_1$, $\al_2$ be embedded 
representatives of $c_1$, $c_2\in H_1(M,\bbZ_2)$ chosen such that
\[
x = \hbox{$\#$ of intersection points of $\al_1$ and $\al_2$} =
c_1 \cdot c_2.
\]
Let $N$ be a regular neighborhood of $\al_1\cup\al_2$ on $M$, with $N$ diffeomorphic to the thrice-punctured sphere
or punctured torus as $x = 0$ or $1$, respectively. Choose an embedded
representative $\beta:S^1\to N$ for $c_1 + c_2$.
Therefore
\[
\begin{array}{rcl}
q_S(c_1+c_2) &=& w(\beta,v) + 1 =
    (w(\al_1,v) + w(\al_2,v) + (x + 1)) + 1\\
&=& (w(\al_1,v) + 1) + (w(\al_2,v) + 1) + x\\
&=& q_S(c_1) + q_S(c_2) + c_1\cdot c_2.
\end{array}
\]

\begin{lemma}
If $\al:S^1\to M$ is an embedded curve on a surface $M$ with spin
structure $S$, and $v_1$, $v_2$ are smooth non-zero vector fields
along $\al$, then $w(\al,v_1) = w(\al, v_2)$ if and only if $v_1$ and $v_2$
alike lift or do not lift along $\al$ to smooth sections of $S$.
\end{lemma}

\begin{pf}
We may assume $M$ is an annulus containing $\al(S^1)$ as the unit circle, 
with spin structure $S_k$ $(k=0 \hbox{ or } 1)$ as in Section
\ref{spin structures}.
Any vector field  $S^1\to\bbC$ is of the form $t\mapsto t^p[f(t)]^2$, where $f$ is smooth and
\[
p = \left\{\begin{array}{ll}
k & \hbox{if $v$ lifts,}\\
1-k &  \hbox{if $v$ does not lift.}
\end{array}\right.
\]
Then, with $w_\al(h_1,h_2)$ defined as the 
winding number (mod 2) of $h_1$ against
$h_2$ (or equivalently, of $h_2/h_1$) along $\al$,
\begin{eqnarray*}
w_\al(v_1,v_2) &=& w_\al(t^p[f_1(t)]^2,t^q[f_2(t)]^2) = w_\al(t^p,t^q)
\equiv p+q \mbox{ (mod 2)}\\
&=& \left\{ \begin{array}{ll}
0 & \hbox{if $v_1$, $v_2$ alike lift or do not lift,}\\
1 & \hbox{otherwise.}
\end{array}\right.
\end{eqnarray*}
But $w_\al(v_1,v_2) = w(\al,v_1) + w(\al,v_2)$, and the result follows.
\end{pf}

\section{Spin structures on hyperelliptic Riemann surfaces}\label{hyper appendix}

Here the spin structures, their corresponding quadratic forms,
and their Arf invariants are computed explicitly for hyperelliptic
Riemann surfaces.

\begin{theorem}\label{hyperelliptic}
Let 
\[
M = \left\{
[x_1,x_2,x_3]\in\bbC{\bbP}^2\,\,
\left|
\,\,x_2^2x_3^{2g-1} = 
\textstyle\prod_{i=1}^{2g+1}(x_1-a_ix_3)
\right.
\right\}
\]
be a hyperelliptic Riemann surface of genus $g$, where
$A = \{\upto{a}{1}{2g+1}\}\subset\bbC$ is a set of $2g+1$ distinct points.
Let $z = x_1/x_3$ and $w = x_2/x_3$.
For each subset $B\subseteq A$, define 
\[
f_B(z) = \prod_{a\in B}(z - a)\andspace \eta_B = f_B(z) dz/w.
\]
Then
\begin{enumerate}
\item
Any differential $\eta_B$ represents
a spin structure in the sense of Theorem~\ref{holo spin theorem}.
\item
The set of $2^{2g}$ meromorphic differentials
\[
\left\{ \eta_B\suchthat B\subseteq A, \# B\le g\right\}
\]
represent the $2^{2g}$ distinct spin structures on $M$.
\item
With $q_B$ the quadratic
form corresponding to $\eta_B$, let $\gam$ be a curve in $M$ whose
projection to the $z$-plane is a Jordan curve which
avoids $\infty$ and $A$, and let
$C\subseteq A$ be the set of branch points which lie in the 
region enclosed by $\gam$ (so $\# C$ is even).  Then
\[
q_B([\gam]) = \# (B\cap C) + \frac{1}{2}\# C\mbox{ (mod 2)}.
\]
\item
With $q_B$ as in (iii),
{\em
\[
\Arf q_B = \left\{\begin{array}{ll}
+1 & \hbox{if\ \ \ $2g - 2\# B + 1 \equiv \pm 1\mbox{ (mod 8)}$},\\
-1 & \hbox{if\ \ \ $2g - 2\# B + 1 \equiv \pm 3\mbox{ (mod 8)}$}.
\end{array}\right.
\]
}
\end{enumerate}
\end{theorem}

\begin{pf*}
{\it Proof of (i)}
Let $P_i = P_{a_i} = [a_i,0,1]$ and 
$P_\infty = [0,1,0]$ be the branch points of
the two-sheeted cover $z:M\to\bbC{\bbP}^1$.
Then the divisor of $\eta_B$ is
\[
 2\left( (g-\# B-1)P_\infty + \sum_{a\in B} P_a\right).
\]
Since this divisor is even, the differential represents a spin structure by
Theorem~\ref{holo spin theorem}.

{\it Proof of (ii).}
Note that there are $\left(2g+1 \atop r\right)$ differentials 
in the $(r+1)^{\hbox{\scriptsize th}}$ row,
totaling $\sum_{r=0}^g\left(2g+1\atop r\right) = 2^{2g}$.  
All but those in the last row are holomorphic.

In order to prove that these differentials represent distinct
spin structures, we first compute the relations on the divisors
of the form
$\sum k_i P_i + k_\infty P_\infty$.
Two such divisors are equivalent if and only if
there is a meromorphic function $M$
whose divisor is their difference.
Since the functions $w$ and $z-a_i$ have respective divisors
\[
\begin{array}{rcl}
(w) &=& P_1 + \cdots + P_{2g+1} - (2g+1)P_\infty,\\
(z-a_i) &=& 2P_i - 2P_\infty,
\end{array}
\]
we have the independent relations
\begin{equation}\label{divisor relations}\begin{array}{rcl}
P_1 + \cdots + P_{2g+1} &\equiv & (2g+1)P_\infty,\\
2 P_i &\equiv & 2 P_\infty \spaceout\spaceout\spaceout {(i=1,\dots,2g+1).}
\end{array}
\end{equation}
To show that there are no other relations independent of these,
let $\sum k_i P_i + k_\infty P_\infty\equiv 0$ be a relation.
Then $\sum k_i = k_\infty$, and by the relations above, we may assume
each $k_i$ is $0$ or $1$.  Hence the general relation may be assumed to
be of the form $D-d P_\infty\equiv 0$, 
where $D$ is a sum of distinct $P_i\in A$, and $d = \# D$.
Let $h$ be a function with divisor $D-d P_\infty$.  Since the only pole of
$h$ is at $P_\infty$, $h$ is a polynomial in $z$ and $w$, so there are 
polynomial functions $f_1$ and $f_2$ of $z$ such that
\[
h(z,w) = f_1(z) + w f_2(z).
\]
Then
\[
2g+1\ge d = -\ord{P_\infty} h = -\ord{P_\infty}(f_1+wf_2)
\ge -\ord{P_\infty} wf_2 = \deg f_2 + 2g + 1.
\]
Thus $d=2g+1$, and $D=P_1+\cdots + P_{2g+1}$,
so no new relation can exist.

We want to show that $\eta_{B_1}$ and $\eta_{B_2}$ represent
identical spin structures if and only if $B_1 = B_2$ or $B_1 = B_2'$, where
the prime notation $C'$ designates the complement $A\setminus C$ in $A$.
If $B_1 = B_2$, then this is clear;  if $B_1 = B_2'$, then
$\eta_{B_2}/\eta_{B_1} = (f_2/w)^2$ 
is a square of a meromorphic function on $M$,
and so $\eta_{B_1}$ and $\eta_{B_2}$ 
represent the same spin structure by Theorem~\ref{holo spin theorem}.

Conversely, suppose that $\eta_{B_1}$ and $\eta_{B_2}$ represent the same 
spin structure.  Then by Theorem~\ref{holo spin theorem},
 $\eta_{B_2}/\eta_{B_1}=h^2$ for some meromorphic function $h$ on $M$.
But 
\[
2(h) = (h^2) =
  (\eta_{B_2}/\eta_{B_1}) = 2((d_2 - d_1)P_\infty + D_2 - D_1),
\]
where $D_1 = \sum_{a\in B_1}P_a$, $D_2 = \sum_{a\in B_2}P_a$, 
$d_1 = \# B_1$, and $d_2 = \# B_2$.
Therefore $(d_2 - d_1)P_\infty + D_2 - D_1\equiv 0$.
By the relations (\ref{divisor relations}), this divisor is
equivalent to
\[
\sum_{a\in B_1\circ B_2}P_a - \# (B_1\circ B_2) P_\infty,
\]
where $B_1\circ B_2$ is the symmetric difference 
$(B_1\cup B_2)\setminus (B_1\cap B_2)$.
Since the relations (\ref{divisor relations}) generate all such relations,
it follows that $B_1\circ B_2$ is either $\emptyset$ or $A$, that is
that $B_1=B_2$ or $B_1=B_2'$.

{\it Proof of (iii).}
It follows from the definition of $q = q_B$ that $q([\gam])$ is the
degree (mod 2) of the map $f(z)/w$ thought of as a map from the curve
$\gam$ on $M$ to $\bbC\setminus\{0\}$.
Let $h = (f/w)^2$.  Then
\[
\deg h = \sum_{h(p) = 0}\ord_{p} h + \sum_{h(p) = \infty}\ord_{p}{h},
\]
the sums being restricted to points within $\gam$. This computes to
\[
\deg h = \#(B\cap C) - \#(B\cup C) = 2(\# (B\cap C) - 
    \textstyle\frac{1}{2}\# C),
\]
which shows that
\[
q([\gam]) = \#(B\cap C) + \textstyle\frac{1}{2}\# C\mbox{ (mod 2)}.
\]

{\it Proof of (iv).}
In order to compute $\Arf q$, we first compute 
$\sum q(\al)$, where $\al$ ranges over $H_1(M,\bbZ_2)$.
Correspondingly, the set of branch points $C$ in the region enclosed
by $\al$ range over
the subsets of $A$ of even cardinality.  Hence $\sum q(\al)$ is the number
of such subsets for which $q(\al)=1$, that is, for which
\[
\#(B\cap C) - \#(B'\cap C) \equiv 2\mbox{ (mod 4)}.
\]
The set of such subsets is
\[
\{R\cup S\suchthat 
  R\subseteq B, S\subseteq B', \# R - \# S \equiv 2\mbox{ (mod 4)}\}.
\]
The cardinality of this set is
\[
\sum q(\al) = \sum_{i-j\equiv 2}
  \left(b\atop i\right)\left(b'\atop j\right),
\]
where $b = \# B$, $b' = \# B'$, and the sum is over $i$ and $j$ with
$i-j\equiv 2\mbox{ (mod 4)}$.

To compute this sum, define
\[
\xi(c,k) = \sum_{i\equiv k}\left(c\atop i\right).
\]
Then 
\begin{eqnarray*}
\sum q(\al) &=& \sum_{i}\left(b\atop i\right)
     \sum_{j\equiv i+2}
\left(b'\atop j\right) = \sum_k\left(b\atop i\right)\xi(b',j+2)\\
&=& \sum_{p=0}^3\sum_{n}\left(b\atop 4n+p\right)\xi(b',p+2) =
\sum_{p=0}^3 \xi(b,p) \xi(b',p+2).
\end{eqnarray*}
Using a fact about Pascal's triangle
\[
\xi(c,k) = 
  2^{(c-2)/2}\left(2^{(c-2)/2} + \cos \textstyle\frac{\pi}{4}(c-2k)\right),
\]
we have
\begin{eqnarray*}\textstyle
\sum q(\al) &=& 2^{(2g-3)/2} \sum_{p=0}^3 
\left( 2^{(b-2)/2} +\cos\textstyle\frac{\pi}{4}(b-2p)\right)
\left( 2^{(b'-2)/2} +\cos\textstyle\frac{\pi}{4}(b'-2p)\right)\\
&=& 2^{g-1}\left( 2^g - \frac{1}{\sqrt{2}}\sum_{p=0}^3
    \cos\textstyle\frac{\pi}{4}(b-2p)  \cos\textstyle\frac{\pi}{4}(b'-2p)
\right)\\
&=& 2^{g-1}\left(2^g - \sqrt{2}\cos\textstyle\frac{\pi}{4}(2g-2b+1)\right)\\
&=&
\left\{\begin{array}{lll}
2^{g-1}(2^g-1) &\hbox{if}& 2g-2b+1\equiv \pm 1\mbox{ (mod 8)},\\
2^{g-1}(2^g+1) &\hbox{if}& 2g-2b+1\equiv \pm 3\mbox{ (mod 8)}.
\end{array}\right.
\end{eqnarray*}
Since $(-1)^t = 1-2t$ for $t = 0$ or $1$,
\[
\Arf q = \frac{1}{2^g}\sum(-1)^{q(\al)} = 
\frac{1}{2^g}(2^{2g} - 2\sum q(\al))
\]
is $+1$ or $-1$ according as $2g-2b+1$ is $\pm1$ or $\pm3$ (mod 8).
\end{pf*}

\section{Group action on spinors}\label{appendix group action}

In this Appendix we outline the proof that $\GL(2,\bbC)$ is the spin
covering group of the linear conformal group (Theorem~\ref{two-fold cover}).

\begin{pf}
Identify $\bbC^3$ with the set $\Gamma$ of
trace-free $2\times 2$ complex matrices via
\[
(x_1,x_2,x_3) \longleftrightarrow\left(
\begin{array}{cc}x_3 & -x_1+ix_2\\ -x_1-ix_2 & -x_3\end{array}
\right)
=X,
\]
and identify ${\bbR}^3\subset\bbC^3$ with
$\Gamma_{\bbR}=\{X\in\Gamma\suchthat X = \trans{\ov{X}}\}.$
The inner product on $\bbC^3$ becomes
\[
X \cdot Y = \textstyle\sum_{1}^{3}x_i y_i =\frac{1}{2}\trace{XY},
\]
and \[
X \cdot X = \frac{1}{2}\trace{X^2} = -\det X,
\] so
$Q\subset\bbC^3$ 
is identified with \[
\Gamma_Q = \{X\in\Gamma\suchthat\det X = 0\}.
\]
Similarly, $\bbC^2$ may be identified with the set $\Delta$ of matrices of the 
form 
\[
\left(\begin{array}{cc}x_1 & x_1\\x_2 & x_2\end{array}\right).
\]
Under these identifications the map $\sigma:\bbC^2\to Q$ becomes
$\sig:\Delta\to\Gamma_Q$ given by $\sig(X) = XJX'$, where 
$J = \left(\begin{array}{rr}0 & -1\\ 1 & 0\end{array}\right)$, and
$X'$ denotes
the classical adjoint
\[
\left(\begin{array}{cc}
a & b\\c & d
\end{array}\right)'
=
\left(\begin{array}{rr}
d & -b\\
-c & a
\end{array}\right)
\]
satisfying $XX' = X'X = (\det X)I$ and $(XY)' = Y'X'$.

In order to satisfy equation (\ref{T commutes}), then $T$ must be defined, for 
$X\in\Gamma$, by
\[
T(A)X = AXA'.
\]
It follows that $T(A)$ is linear and maps $\Gamma$ to itself, and that
$T:\GL(2,\bbC)\to\GL(3,\bbC)$ is a homomorphism 
with kernel $\{\pm I\}$.
That $T$ restricts as indicated follows from the equation
\[
T(A)X \cdot T(A)Y = (\det A)^2 X \cdot Y
\]
and the fact that $T(A)(\Gamma_{\bbR}) = \Gamma_{\bbR}$ for 
$A\in{\bbR}^*\times\SU(2)$.
\end{pf}

\section{The pfaffian}\label{pfaffian}
Here we recall some basic facts about skew-symmetric forms.

\begin{definitionU}
A bilinear form $A$ on a vector space $V$ of dimension $n$ is 
{\em skew-symmetric} if
\[
 A(v_1,v_2) + A(v_2,v_1) = 0\mbox{    for all $v_1$, $v_2\in V$},
\]
or alternatively, if the matrix $A$ for $A$ satisfies
\[
A+\trans{A} = 0.
\]
\end{definitionU}
The space of skew-symmetric bilinear forms is $\bigwedge^2(V^*)$.
The pfaffian is a function on skew-symmetric forms whose square
is the determinant.

\begin{definitionU}
For $A\in\bigwedge^2(V^*)$, the {\em pfaffian} of $A$ is
{\em
\[
\pfaffian A = \left\{
\begin{array}{cl}
\frac{1}{m!}\overbrace{(A\wedge\dots\wedge A)}^{m\rm\;times} &
    \mbox{if $\dim(V)=2m$ is even,}\\
0 & \mbox{if $\dim(V)$ is odd.}
\end{array}
\right.
\]
}
\end{definitionU}

For a matrix $(a_{ij})$ of $A\in\bigwedge^2(V^*)$
in the basis $\{\upto{e}{1}{m}\}$
the pfaffians for $m=2$, $m=4$, and $m=6$ are respectively
\[
\begin{array}{c}
a_{12},\\
a_{12}a_{34} - a_{13}a_{24} + a_{14}a_{23},\\
a_{12} a_{34} a_{56} - a_{12} a_{35} a_{46} + a_{12} a_{36} a_{45}-
a_{13} a_{24} a_{56} + a_{13} a_{25} a_{46} -\\ a_{13} a_{26} a_{45}+
a_{14} a_{23} a_{56} - a_{14} a_{25} a_{36} + a_{14} a_{26} a_{35}-
a_{15} a_{23} a_{46} +\\ a_{15} a_{24} a_{36} - a_{15} a_{26} a_{34}+
a_{16} a_{23} a_{45} - a_{16} a_{24} a_{35} + a_{16} a_{25} a_{34}.
\end{array}
\]
The general pfaffian of a $2m\times 2m$ matrix has
$(2m)!/(2m!) = 1\cdot 3\cdot 5\cdot \cdots \cdot(2m-1)$ terms.

\begin{lemmaU}
The rank of a skew-symmetric matrix is even.
\end{lemmaU}

\begin{pf}
Let $A$ be an $m \times m$ skew-symmetric matrix with rank $r$.  The
proof is by induction on $m$.  In the case $m=1$, then $A = (0)$ with
even rank $0$.  Assume for some $n$ that the lemma is true for all
skew-symmetric matrices smaller than $A$.  If $n$ is odd, then
\[
\det A = \det \trans{A} = \det(-A) = (-1)^n \det A = - \det A,
\]
so $\det A=0$ and $A$ has a non-zero kernel.  If $n$ is even,
then $A$ also has  a non-zero kernel unless it has full --- hence even
--- rank $r=n$.  So in either case we may assume $A$ has a non-zero kernel.

Let $\upto{v}{1}{n-r}$ be a basis for $\ker A$, and let 
$\upto{v}{1}{n-r}$, $\upto{w}{1}{r}$ be an extension of this basis
to a basis for $\bbC^n$.  Let $P$ be the
$n \times n$ matrix with these vectors as columns.  Then $\trans{P}AP$ is of
the form
\[
\trans{P}AP = \left(
\begin{array}{c|c}
0 & 0 \\ \hline
0 & A_0
\end{array}
\right),
\]
where $A_0$ is an $r \times r$ matrix of rank $r<n$.  Moreover, 
\[
\trans{(\trans{P}AP)} = \trans{P}\trans{A}P = -(\trans{P}AP),
\]
so $\trans{P}AP$, and hence $A_0$ is skew-symmetric.  By the induction
hypothesis, $r=\rank A$ is even, since it is the rank of the
smaller skew-symmetric matrix $A_0$.
\end{pf}

\section{Elliptic functions}\label{elliptic}
For reference, here are some standard
notations and facts about elliptic functions used in this paper
(see for example \cite{DuVal}, \cite{Erdelyi}).

\skipline
{\em Lattices.}  A non-degenerate lattice $\Lambda$ is {\em real} if
$\Lambda = \ov{\Lambda}$.  There are two kinds of real lattices:
\begin{enumerate}
\item rectangular:  generators $\om_1\in{\bbR}$ and $\om_3\in i{\bbR}$ can 
be chosen for $\Lambda$.
\item rhombic:  generators $\om_1$ and $\om_3=\ov{\om}_1$ can be chosen
for $\Lambda$.
\end{enumerate}
For any lattice with generators $\om_1$, $\om_3$, let $\om_2 = -\om_1-\om_3$.

\skipline
{\em The Weierstrass $\wp$ function:}  Given a lattice $\Lambda$
generated by $\om_1$ and $\om_3$, the elliptic function $\wp$
on $\bbC/\Lambda$ satisfies the differential equation
\[
  (\wp')^2 = 4 \wp^3 - g_2 \wp - g_3 = 4(\wp-e_1)(\wp-e_2)(\wp-e_3),
\]
where
\[
\begin{array}{l}
e_i = \wp(\om_i) \spaceout (i=1,2,3),\\
e_1+e_2+e_3=0,\\
g_2 = -4(e_1 e_2+e_1 e_3+e_2 e_3),\\
g_3 = 4 e_1 e_2 e_3.
\end{array}
\]
The function $\wp$ has a double pole at $0$ and two simple zeros which
come together only on the square lattice; $\wp'$ has a triple pole at
$0$ and three simple poles at $\om_1$, $\om_2$, $\om_3$.

The function $\wp$ is even; $\wp'$ is odd.  On a horizontal rectangular
lattice, $\wp(\ov{u}) = \ov{\wp(u)}$; on a horizontal 
square lattice, $\wp(iu) = -\wp(u)$.

The expansion for $\wp$ at $0$ is
\[
\wp(u) = \frac{1}{u^2} + \frac{g_2}{20} u^2 + \dots.
\]

A useful property of $\wp$ is the following special case of the
addition formula ($\{i,j,k\}$ is any permutation of $\{1,2,3\}$):
\begin{equation}\label{addp}
\wp(u\pm\om_i) = e_i + \frac{(e_i-e_j)(e_i-e_k)}{\wp(u)-e_i}.
\end{equation}

\skipline
{\em The Weierstrass $\zeta$ function:}  The $\zeta$ function is
defined by
\[
\zeta(u) = -\int \wp(u)du,
\]
with the constant of integration chosen so that 
$\lim_{u\rightarrow 0}\zeta(u)-u^{-1} = 0$.
With $\eta_i = \zeta(\om_i)$ $(i=1,2,3)$, properties of $\zeta$ include:
\[
\begin{array}{l}
\eta_1+\eta_2+\eta_3=0,\\
\zeta(u+2\om_i) = \zeta(u) + 2 \eta_i\spaceout (i=1,2,3),\\
\zeta \mbox{ is an odd function.}
\end{array}
\]
Legendre's relation is that

\begin{equation}\label{Legendre}
\eta_1 \om_3 - \eta_3 \om_1 = i\pi/2.
\end{equation}
A form of the quasi-addition formula for $\zeta$ is
\begin{equation}\label{zeta}
\zeta(u-v) - \zeta(u)+\zeta(v) = \frac{1}{2}
     \left(\frac{\wp'(u)+\wp'(v)}{\wp(u)-\wp(v)}\right).
\end{equation}

\skipline
A useful property of elliptic functions which can also be stated in more
generality is the following:
Let $f$ be an elliptic function with poles of order at most $2$, with
no residues, and with principal parts
\[
\frac{a_1}{(u-\al_1)^2},\dots,\frac{a_n}{(u-\al_n)^2}.
\]
Then
 
\begin{equation}\label{standard}
f(u) = b+\sum a_i\wp(u-a_i)
\end{equation}
for some $b$, because the difference $f(u)-\sum\al_i\wp(u-\al_i)$ 
has no poles and hence is constant.

\section{Klein bottles: conformal type, spin structure, periods}
\label{klein}

Here we show that
the torus covering a Klein bottle must 
have the conformal type of the complex plane modulo a rectangular lattice,
we compute the order-two deck transformation for the covering, and we show that the spin structure
on such a torus must be untwisted. (This can also be seen from purely topological 
considerations.)

\begin{theorem}\label{kleintype}
Let $X:K'\longrightarrow\bbR^3$ be a complete minimal 
immersion of a punctured Klein bottle with finite total curvature,
$\pi:T\longrightarrow K=\ov{K'}$ the oriented two-sheeted
covering by a torus $T$, and $I:T\longrightarrow T$ the order-two
orientation-reversing deck transformation for this cover.
Then we have the following.
\begin{enumerate}
\item $T$ is conformally equivalent to $\bbC/\Lambda$, where
$\Lambda$ is a rectangular lattice with generators
$2\om_1\in\bbR$ and $2\om_3\in i\bbR$.
\item On this torus, the deck transformation $I$ may be chosen to be 
$I(u)=\bar{u}+\om_1$.
\item With this choice, the admissible spin structures 
are those represented by $(\wp(u)-\wp(\om_2))du$ and
$(\wp(u)-\wp(\om_3))du$.
\item If $(s_1,s_2)$ is the spinor representation of $X\circ\pi$ on $T$,
the period conditions reduce to the conditions 
$\int_{\gam_1} s_1^2 = 0$ and $\int_{\gam_1} s_1s_2 = 0$ 
along a closed curve $\gam_1$ parallel to $\om_1$.
\end{enumerate}
\end{theorem}

\begin{pf*}
{\it Proof of (i) and (ii)}  Let $\Lambda_0$ be a lattice such that
$T=\bbC/\Lambda_0$.  Since every conformal map from $T$ to $T$ must be
linear in the standard coordinate $u$ on $\bbC$ 
and since $I$ is anti-conformal, $I(u) = \al\bar{u}+\beta$ for 
some $\al$, $\beta\in\bbC$.  The periodicity of $I$ and $I^{-1}$ implies
that $\al\ov{\Lambda_0}\subseteq\Lambda_0$ and 
$\ov{\al}^{-1}\ov{\Lambda_0}\subseteq\Lambda_0$.  These together imply that
$\al\ov{\Lambda_0} = \Lambda_0$.  Choose $\gam\in\bbC$ satisfying
$|\gam|=1$ and $\ov{\gam}/\gam=\al$;  the rotated lattice 
$\Lambda = \gam\Lambda_0$ satisfies $\ov{\Lambda} = \Lambda$
(a so-called {\em real} lattice).
Hence $\Lambda$ is either rectangular with generators
$2\om_1\in\bbR$, $2\om_3\in i\bbR$, or $\Lambda$ is rhombic with
generators $2\om_1$ and $2\om_3=2\ov{\om}_1$.
On $\bbC/\Lambda$ we have
$I(u) = \al \bar{u}+\beta$ for some new $\al$, $\beta\in\bbC$.
As before, $\al\ov{\Lambda} = \Lambda$, but $\ov{\Lambda}=\Lambda$,
so $\al=\pm 1$.  If $\al=-1$, replacing $\Lambda$ by $i\Lambda$
preserves its reality, and changes $\al$ to $1$.

With $\al=1$, the condition that $I$ is involutive is that
$\beta + \ov{\beta}\in\Lambda$. 
By the change of coordinate $u\mapsto u - i\;\im \beta$, it can be
assumed that $\beta\in\bbR$.  Then the involutive condition
is that $2\beta\in\Lambda$.  
If $\beta\in\Lambda$ then $0$ is a fixed point of $I$.  Hence
$\beta \equiv \om_1$ (rectangle) or $\beta = \om_1+\om_3$ (rhombus).
In the latter case, $\om_1$ is a fixed point of $I$, 
so the only admissible case is the rectangle, with $I(u) = \ov{u}+\om_1$.

{\it Proof of (iii).} The compatibility condition in
Theorem~\ref{nonorientable theorem}
demands that $I^*I^*(s)=-s$ for any section $s$ of the spin structure.
A computation shows that this condition is met only 
for the two spin structures named.

{\it Proof of (iv).}
Let $\gam_1$ and $\gam_3$ be respectively the closed curves 
$t\mapsto\om_1t/|\om_1|+c_1$ and $t\mapsto\om_3 t/|\om_3|+c_2$,
$(0\le t\le 2)$, where
$c_1$, $c_2\in\bbC$ are chosen so that the curves do not pass through
any ends.  Then $I(\gam_1) = \gam_1$, $I(\gam_3) = -\gam_3$.
The period conditions are
\[
\int_{\gam_k} s_1^2 = \ov{ \int_{\gam_k} s_2^2}
\andspace
\int_{\gam_k} s_1 s_2 \in i\bbR\spaceout (k=1,3).
\]
With $I$ as above, under the double-cover assumption 
\[
(s_1,s_2) = \pm(i \ov{I^*s_2}, -i \ov{I^*s_1}),
\]
we have
\[
\int_{\gam_3} s_1^2 = \int_{\gam_3} -\ov{I^*s_2^2} =
  -\ov{\int_{I(\gam_3)} s_2^2} = \ov{\int_{\gam_1} s_2^2}
\]
\[
\int_{\gam_3} s_1s_2 = \int_{\gam_3} \ov{I^*s_1s_2} =
  \ov{\int_{I(\gam_3)} s_1s_2} = -\ov{\int_{\gam_3} s_1s_2},
\]
so the period conditions are automatically satisfied for $k=3$.
Moreover, we also have
\[
\int_{\gam_1} s_1^2 = \int_{\gam_1} -\ov{I^*s_2^2} =
  -\ov{\int_{I(\gam_1)} s_2^2} = -\ov{\int_{\gam_1} s_2^2}
\]
\[
\int_{\gam_1} s_1s_2 = \int_{\gam_1} \ov{I^*s_1s_2} =
  \ov{\int_{I(\gam_1)} s_1s_2} = \ov{\int_{\gam_1} s_1s_2}
\]
and the first two period conditions (\ref{period1}) become
\[
\int_{\gam_1} s_1^2 = 0 \andspace \int_{\gam_1} s_1s_2 = 0
\]
(this amounts to three real conditions because, under the above assumption,
the second integral is automatically real).
\end{pf*}

%
%

%
%

\typeout{_______________________________________________ bib.tex}

\frenchspacing  

\end{document}